\documentclass[aos]{imsart}

\usepackage{amsfonts, amsmath, amsthm, amssymb, textcomp,enumerate,
  array, latexsym, paralist, mathtools,bbm, color}
\RequirePackage[numbers]{natbib}
\RequirePackage[colorlinks,citecolor=blue,urlcolor=blue]{hyperref}%% uncomment this for coloring bibliography citations and linked URLs
\RequirePackage{graphicx}%% uncomment this for including figures
\usepackage{paralist}
\usepackage{subfigure}
\usepackage{float}
\startlocaldefs
\startlocaldefs
\theoremstyle{plain}

\newtheorem{theorem}{Theorem}[section]
\newtheorem{lemma}[theorem]{Lemma}
\newtheorem{remark}[theorem]{Remark}
\newtheorem{corollary}{corollary}
\newcommand{\mynorm}[1]{ \left\| #1 \right\| }
\theoremstyle{remark}

\newcommand{\R}{\mathbb{R}}
\newcommand{\E}{\mathbb{E}}

\newcommand{\skewF}{\overline{F}}
\newcommand{\skewQ}{\overline{Q}}
\newcommand{\geomQ}{Q_G}
\newcommand{\hgeomQ}{\widehat{Q}_G}

\endlocaldefs

\begin{document}

\begin{frontmatter}
\title{On the use of the M-quantiles for outlier detection in multivariate data}
%\title{A sample article title with some additional note\thanksref{T1}}
\runtitle{M-quantiles and outlier detection}
%\thankstext{T1}{A sample of additional note to the title.}

\begin{aug}
%%%%%%%%%%%%%%%%%%%%%%%%%%%%%%%%%%%%%%%%%%%%%%%
%% Only one address is permitted per author. %%
%% Only division, organization and e-mail is %%
%% included in the address.                  %%
%% Additional information can be included in %%
\author{\fnms{Sajal}~\snm{ Chakroborty}\ead[label=e1]{sajal.chakroborty@ttu.edu}}
\and
\author{\fnms{Ram}~\snm{Iyer}\ead[label=e2]{ram.iyer@ttu.edu}\orcid{0000-0001-8908-0892}}
\and
\author{\fnms{A.~Alexandre}~\snm{Trindade}\ead[label=e3]{alex.trindade@ttu.edu}}
%%%%%%%%%%%%%%%%%%%%%%%%%%%%%%%%%%%%%%%%%%%%%%
%% Addresses                                %%
%%%%%%%%%%%%%%%%%%%%%%%%%%%%%%%%%%%%%%%%%%%%%%
\address{Department of Mathematics and Statistics,
	Texas Tech University, USA \printead[presep={,\ }]{e1,e2,e3}}
	
%	\address[B]{Department of Mathematics and Statistics,
	%		             Texas Tech University, USA\printead[presep={,\ }]{e2,e3}}
\end{aug}

\begin{abstract}
Defining a successful notion of a multivariate quantile has been an open problem for more
than half a century, motivating a plethora of possible solutions. Of
these, the approach of \cite{chakraborty1996transformation} and
\cite{koltchinskii1997m} leading to M-quantiles, is very appealing for
its mathematical elegance -- combining elements of convex analysis and
probability theory. The key idea is the description of a convex
function (the K-function) whose gradient (the K-transform) is in 
one-to-one correspondence between all of $\R^d$ and the unit ball in
$\R^d$. By analogy with the $d=1$ case where the K-transform
is a cumulative distribution function-like object (an M-distribution), the fact that its
inverse is guaranteed to exist lends itself
naturally to providing the basis for the definition of a quantile
function for all $d\geq 1$. Over the past twenty years the resulting M-quantiles
have seen applications in a variety of fields, primarily for the
purpose of detecting outliers in multidimensional spaces.
In this article we prove that for odd $d\geq 3$, it is not the gradient but a
poly-Laplacian of the K-function that is (almost
everywhere) proportional to the density
function. For $d$ even one cannot establish a differential equation connecting
the K-function with the density.
These results show that usage of the K-transform for outlier
detection in higher odd-dimensions is in principle flawed,
as the K-transform does not originate from inversion of a true M-distribution.
We demonstrate these conclusions in two dimensions
through examples from non-standard asymmetric distributions. Our
examples illustrate a feature of the K-transform whereby regions in
the domain with higher density map to larger volumes in the co-domain,
thereby producing a magnification effect that moves inliers closer to
the boundary of the co-domain than outliers. This feature obviously
disrupts any outlier detection mechanism that relies on the inverse K-transform.
\end{abstract}


\end{frontmatter}
%%%%%%%%%%%%%%%%%%%%%%%%%%%%%%%%%%%%%%%%%%%%%%
%% Please use \tableofcontents for articles %%
%% with 50 pages and more                   %%
%%%%%%%%%%%%%%%%%%%%%%%%%%%%%%%%%%%%%%%%%%%%%%
\tableofcontents
%%%%%%%%%%%%%%%%%%%%%%%%%%%%%%%%%%%%%%%%%%%%%%
%%%% Main text entry area:
\newpage
\section{Introduction}\label{intro} 
The subject of outlier detection has a long history in the statistical
and data-mining literature. This harks at the difficulty of the
enterprise, which reflects the fact that a precise definition is
ellusive. Thus, while there is general agreement that an outlier is an observation that
deviates so much from the others as to arouse suspicion,
or appears to be inconsistent with the remainder of the data, there is
much debate on how to measure the actual ``deviation'' or
``inconsistency''~\citep{ben2005outlier}. 

For univariate data that can plausibly be modeled as a realization of
a random sample (independent and identically distributed, or i.i.d.), most detection methods
rely on the quantile function of the data
generating mechanism. In the taxonomy of the field one can distinguish
between parametric methods which classify outliers as those lying on
the tails of the distribution, and nonparametric methods that rely on inter-quantile
ranges~\citep{tukey1977exploratory}. For multivariate data, the inherent difficulty in
defining quantiles in multidimensional spaces has shifted the focus
to measuring the distance of individual points from the center of the
data cloud. Thus one finds a substantial body of literature on such
\emph{distance functions}, or their inverses, \emph{depth functions};
see for example \cite{Mosler2013} for a survey, and
\cite{Serfling2000depth} for formal definitions and theoretical development. 

Quantiles are calculated by inverting the
cumulative distribution function (cdf).
Mathematically, for any random variable $X$ defined on $\R$, the most generally accepted definition of the quantile function
at probability level $s$, is the mapping $Q:[0,1]\mapsto\R$, given by 
%%%%%%% Equation 1 %%%%%%%%%%%%%%%%%%%%%%
\begin{equation}\label{eqn: QNTL 1.1}
	Q(s) = \inf \{x: F(x) \geq s\},
\end{equation}
%%%%%%%%%%%%%%%%%%%%%%%%%%%%%%%%%%%%%%%%%
where $F(x)=P(X\leq x)$ denotes the cdf of $X$, corresponding to a
left-to-right accumulation of mass. In one dimension, the only other
essentially different mass accumulation would be center-outward (with
``center'' taken to be, e.g., the median). We term $\skewF(x) = 2\,F(x) -1$ the
\emph{signed} cdf. It is the so-called \emph{M-distribution} arrived at by
\cite{koltchinskii1997m} when attempting to extend the notion of a quantile. 
\cite{li2018edf} notes that goodness-of-fit
tests like Kolmogorov-Smirnov based on mass accumulation dictated by
$|\skewF(x)|$, which corresponds to a center-outward ordering, are more
powerful at detecting scale differences. Such a center-outward
ordering also features prominently in \cite{hallin-etal-2021}, who
define a multivariate cdf (and associated quantile function) based on
the intriguing concept of transportation-based ranks and signs.
One may therefore define an
alternate quantile function based on the mapping  $\skewQ:[-1,1]\mapsto\R$, given by 
\begin{equation} \label{eq:symmetric cdf}
  \skewQ(s) =  \inf \{x: \skewF(x) \geq s \}.
\end{equation}

As there is no natural ordering in multidimensional spaces (dimension
$d$ in general), it is
difficult to extend the notion of a quantile. 
%Geometric Median (GM) or Torricelli point ( \cite{cieslik2006shortest}) could  to define multivariate quantiles. 
Early contributions to this field of study were made by
\cite{haldane1948note}, \cite{brown1983statistical},
\cite{oja1983descriptive}, \cite{breckling1988m}, and \cite{chaudhuri1992multivariate}, among
others. \cite{serfling2002quantile} provides a comprehensive overview and an
exhaustive comparison of the different methodologies. Although it is
not our intent to give an updated survey here, we
highlight two approaches that serve to give some idea of the
difficulties involved.    \cite{breckling1988m} define an M-quantile through the notion of an
M-estimator and its associated loss and influence functions. Although
this leads to sensible definitions for the median, e.g., under
radially symmetric loss we obtain the geometric median (discussed in
detail below), the quantile extension
is cumbersome and suffers from a high computational burden. Another
definition for finite as well as for infinite-dimensional spaces
(e.g., functional data) was
proposed by \cite{fraiman2012quantiles}. Their approach however, which
builds on earlier 
ideas by \cite{kong2012quantile} and others that rely on directional
projections of the probability measure onto the unit sphere,  leads to quantiles
that are not affinely equivariant. (Instead of plotting quantile
contours, they proposed computing principal quantile directions in
order to detect outliers.)

A common thread in these various attempts to generalize the
notion of a quantile, is the imposition of the requirement that the
specialization to the univariate case corresponds to the quantile
function $Q$ or $\skewQ$. Any such attempt must also by necessity
be grounded on the simpler problem of leading to the special case of a
well-defined median. In this regard, it has long been noted
that  $Q(0.5)$ minimizes the mean absolute error (\cite{stroock2011}), 
\[
  Q(0.5) = \underset{c \in \mathbb{R}}{\arg \min}\, \E |X-c| = \underset{c \in \mathbb{R}}{\arg \min}\, \E (|X-c| - |X|),
\]
where the last expression is to be used in  cases when the
expectation does not exist. By extension, the minimizer of the multivariate mean absolute error is typically defined as the multivariate median:
\begin{equation} \label{eqn: QNTL 1.2}
  \text{Med}_{p} = \arg\min_{s\in\R^d}\E (\|X-s\|_p - \|X\|_p),
\end{equation}
where $\|\cdot\|_p$ denotes the $L^p$-norm in
$\mathbb{R}^d$. In the literature, one finds both the $L^2$-norm
(Euclidean distance) and
$L^1$-norm (Manhattan distance) used in \eqref{eqn: QNTL 1.2}; \cite{Oja2013} calls the former the \emph{spatial
  median}, and the latter the \emph{marginal
  median} (vector of marginal medians). The marginal median, first studied by \cite{puri1971} and later
investigated for outlier detection in
\cite{chakraborty2001affine}, suffers from lack of invariance with respect to coordinate rotations, or more
precisely, actions of the special orthogonal group $SO(d)$. Other names in common usage when the 2-norm is employed in
\eqref{eqn: QNTL 1.2}, include \emph{geometric median} and \emph{$L^1$ median}
(\cite{small1990,lopuhaa1991,rizzi1998}). It is easily shown that
geometric median is invariant with respect to coordinate rotations and
translations of the origin. Henceforth in this article, we will
consider the $L^2$-norm in Equation \eqref{eqn: QNTL 1.2} as we attempt to 
generalize and clarify various notions.

%\cb
%[THIS STUFF CAN BE TAKEN OUT BECUASE IT'S WELL-KNOWN, see Oja (2013)] There is a danger to usi%ng the 1-norm (as done in \cite{chakraborty2001affine}), then one can easily prove that:
%\begin{lemma}
%	Assuming $E \|X\|_1 < \infty$, consider the multivariate one-norm median: $\overline{x}_{GM}^1 = \underset{s \in \mathbb{R}^d}{\arg \min} \,E \|X-s\|_1.$ Then,
%	\[
%	\overline{x}_{GM}^1 = (Q_{X_1}(0.5),\,Q_{X_2}(0.5),\,\cdots,\,Q_{X_d}(0.5)). 
%	\] In other words, the multivariate one-norm median is simply the marginal median. 
%      \end{lemma}
%\nc

Chaudhuri \cite{chaudhuri1996geometric} extended the notion of the
geometric median and defined a quantity that he referred to as the
\emph{geometric quantile}. For $v\in B_1(0)$, the open unit ball in
$\mathbb{R}^d$, this is given by 
%%%%%%%%%% Equation 03 %%%%%%%%%%%%%%%%%
\begin{equation}\label{eqn: QNTL 1.3}
 \geomQ(v) = \underset{s \in \mathbb{R}^d}{\arg \min}\; \E \{\|X-s\|_2 + \langle v,X-s\rangle -\|X\|_2 - \langle v,X\rangle \},
\end{equation}
%%%%%%%%%%%%%%%%%%%%%%%%%%%%%%%%%%%%%%%
where $\langle\cdot,\cdot\rangle$ denotes the usual Euclidean inner product. The justification given for
referring to this as a \textit{quantile} is that
$\geomQ(0)$ corresponds to the  geometric
median, i.e., \eqref{eqn: QNTL 1.2} with $p=2$. Contemporaneously,
Koltchinskii \cite{koltchinskii1997m} developed a function
(henceforth \emph{K-function}), which in principle resembles the integral of a multivariate
cdf,
\begin{equation} \label{eq:Kolchinski function}
  f(s) = \E[\|s-X\|_2-\|X\|_2].
\end{equation}
If we now consider the subdifferential function $\partial f(s)$, which we call the
\emph{K-transform} (defined in \eqref{eqn: Koltchinskii transform} with explicit reference to the underlying measure
$\mu$),  Koltchinskii's \emph{M-quantile} is then
the inverse K-transform,  which constitutes
his attempt at defining a multivariate quantile. 
To relate Chaudhuri's notion of a quantile to Koltchinskii's, note that: 
\begin{equation}\label{eqn:geom-quantile}
 \geomQ(v) = \underset{s \in \mathbb{R}^d}{\arg \min}\; \E \{\|X-s\|_2 + \langle v,X-s\rangle \} =  \underset{s \in \mathbb{R}^d}{\arg \max}\; \E \{ \langle
             v,s-X\rangle - \|s-X\|_2 \}.
\end{equation}
Now, since $f(s)$ is convex in $s$, we see that the function
$\geomQ(v)$ is simply the Fenchel
conjugate of $f(s)$.

The noteworthy difference between these  two approaches is that,
whereas \cite{chaudhuri1996geometric} directly defines
$\geomQ(v)$ as a quantile without relating it
to the inverse of a suitably defined cdf, \cite{koltchinskii1997m}
starts with $f(s)$, defines  $\partial f(s)$ as an \emph{M-distribution}
(which coincides with the signed cdf $\skewF(s)$ in the unidimensional case), and uses results from convex analysis to obtain $\geomQ(v) = \partial f^{-1}(v)$.
An immediate corollary of this fact is that the geometric quantile
uniquely determines the underlying distribution
(\cite{koltchinskii1997m}). \cite{chakraborty2001affine} claims that
this property may be used for outlier detection, in that values of
$\|\geomQ(v)\|_2$ closer to one are indicative of
outliers. However, \cite{chen2008} notes that this claim may not hold
in $\R^2$, by demonstrating that outliers may have a smaller value of $\|\geomQ(v)\|$ than inliers for non-standard distributions such as mixtures of Gaussians. However, $d = 2$ is a special case for the K-transform as may be seen in \cite[Theorem 3.2]{chakraborty2001affine}, and it is unclear whether outlier detection may be performed for higher dimensional data. 

The concept of geometric quantiles has been a very active area of research. Although it is
not very useful for outlier detection beyond two dimensions as we show in this article, it is
aesthetically pleasing due to connections with convex analysis. The mathematical elegance
also permits generalizations to infinite dimensional Banach spaces. With equation \eqref{eqn: QNTL 1.3} as the population geometric
quantile, the empirical geometric quantile for a $d$-dimensional
sample of $n$ data points $\{x_1,\ldots,x_n\}$ is
%%%%%%%%%%% Equation 03 %%%%%%%%%%%%%%
\begin{equation}\label{eqn: EQNTL}
\hgeomQ(v) = \underset{s \in \mathbb{R}^d}{\arg \min}\; \frac{1}{n}\,\sum_{i=1}^n \{\|x_i-s\|_2 + \langle v,x_i-s\rangle -\|x_i\|_2 - \langle v,x_i\rangle \}.
\end{equation}
\cite{chaudhuri1996geometric} and \cite{chakraborty2001affine} have derived a Bahadur-type representation for $\hgeomQ(v)$, and proved Glivenko-Cantelli and Donsker-type results. \cite{chakraborty2014spatial} have extended these results for data on infinite dimensional Banach spaces.
Although the quantiles in \eqref{eqn: QNTL 1.3}  are equivariant under
rotations of the data cloud, they are not equivariant under affine
transformations. \cite{chakraborty1996transformation} develop the
notion of an affine equivariant multivariate median for finite
dimensional data, and \cite{chakraborty2001affine} extend this notion
to $L^p$ spaces. Computation of $\geomQ(v)$ has however been a
stumbling block beyond the simplest of distributions. In
Appendix~\ref{sec:app-geomQ-algo} we present an algorithm for this purpose.

Several extensions and applications of these M-quantiles have appeared in
the literature. Their main use has been to effect outlier detection in
areas as varied as manufacturing processes, image processing, economic models, risk
assessment, and biostatistics. With very few exceptions, these are all closely
related to the geometric quantiles of \cite{chaudhuri1996geometric} or
the inverse K-transform of \cite{koltchinskii1997m} (which are
essentialy identical as explained earlier). They were proposed as a means to detect outliers in
two-dimensional data in \cite{chakraborty2001affine} and
\cite{chaouch2010design}. However, the underlying distributions were
only chosen to be the Gaussian, and
hence this paradigm may not generalize well. It has also been proposed
that regions of low probability density should be associated with outliers in data (\cite{chandola2009,aggarwal2015}). 

Various tools called projection quantiles (\cite{mukhopadhyay2011}), quantile contour plots (\cite{chakraborty2001affine,chaouch2010design,girard2017intriguing}), multivariate
quantile-quantile (QQ) plots
(\cite{chakraborty2001affine,marden2004positions,Dhar-etal-2014}), and
spatial depth functions (\cite{serfling2002depth}), have been
developed based on these M-quantiles.
Some applications of these tools include \cite{chen2008} who use them
for outlier detection in images, \cite{church2008spatial} for image
denoising through spatial depth, and \cite{tarabelloni2018statistical}
who devised a
calibration methodology to improve the quality of the
electrocardiograms. The concept of M-quantiles has also found uses in
risk assessment (\cite{herrmann2018multivariate}), control charts in manufacturing processes
(\cite{li2013nonparametric}), climate modeling
(\cite{fieldgenton2006gh}), and development of economic models
(\cite{aryal2016identifying}).

In this article, we will investigate the connection between the probability density
          function (pdf) and the K-function in order to determine the
          utility of the K-transform for outlier detection in
          multivariate data. For $d = 1$, note that the
          K-transform $\nabla f(s)=f'(x)=\skewF(x)$, which 
          reduces to the signed cdf. (This was the original motivation of
          \cite{koltchinskii1997m} in defining $\partial f(s)$ as an
          M-distribution.) Thus $f''(s) = 2\rho(s)$ almost everywhere, where $\rho(\cdot)$ is the underlying pdf of $X$. As we will show, this is the only
          case where differentials of the K-transform lead directly to
          the pdf. Specifically, and with regard to continuous
          measures only for which the subdifferential reduces to the
          ordinary differential (or gradient), $\partial f(s)=\nabla f(s)$, this article makes the following contributions.
\begin{itemize}
	\item For odd dimensions $d \geq 3$, we show it is in fact the
          poly-Laplacian of order $(d+1)/2$ of $f(s)$ that is (almost everywhere) equal to
          a scalar multiple of $\rho(s)$. If one is to believe
          that regions of low probability density are indicative
          of outliers and anomalies in data (e.g.,
          \cite{chandola2009,aggarwal2015}), then this 
          shows that $\nabla f(s)$ does not directly lead to 
          the discovery of regions of low or high probability density in
          high odd dimensional spaces. Thus the M-distribution of \cite{koltchinskii1997m}
          cannot be considered as a distribution in any sense in
          higher dimensional spaces.

	\item For dimension $d=2$ we elucidate the reason for the observation of \cite{chen2008} mentioned above, that values of
$\|\geomQ(v)\|_2\approx 1$  are not necessarily indicative of
outliers. It was shown by \cite{koltchinskii1997m} that the
K-transform is in a one-to-one correspondence between $\R^d$ and the
open unit ball in $\R^d$. However, as we will show, the K-transform
has a zoom-in effect whereby regions of high density occupy a higher
volume in the co-domain, and conversely, regions of low density occupy
a smaller volume. Due to this effect, the distance of a point from the
origin in the co-domain does not necessarily imply that the point lies
in a low density region of the distribution (as would an outlier).   

        \item	Exploiting the connections with convexity theory, we
          propose an algorithm to compute the geometric
          quantile $\geomQ(v)$. The algorithm, based on the MM
          principle \citep{lange2016mm}, allows one to easily
          verify the numerical results obtained by
          \cite{chakraborty2001affine} and \cite{chen2008}, and to
          compute the inverse of the K-transform, even for highly
          nonstandard distributions, and in any finite dimensional vector space. 
\end{itemize}

%\noindent \textbf{Organization of the paper}
These results lead to the main conclusion of this article: that the K-transform cannot
theoretically be considered to be a distribution for dimensions $d
\geq 2$, and consequently, that its inverse cannot be viewed as a quantile function. 
The remainder of the paper is structured as follows. The K-transform
and its inverse are defined and discussed in Section
\ref{sec:K-transform}. In Section \ref{sec:results}, we present formal
results of the points made above concerning the K-transform and empirical geometric
          quantile, and discuss their significance. In  Section
          \ref{sec:zoom-effect}, we present numerical examples to
          illustrate the zoom-in effect of the K-transform, and
          comment on its potential use for outlier detection in two
          dimensional data. The proofs of all formal results are
          relegated to the appendix. 

Before we proceed to the formal definitions, we note that a related approach of using a convex function to define center-outward distribution and quantiles has seen attention in the literature lately \cite{ chernozhukov2017monge,hallin-etal-2021,ghosal2022}. This approach uses the celebrated Brenier-McCann theorem that asserts the existence of a convex function whose gradient transports a reference distribution absolutely continuous with respect to the Lebesgue measure on  $\mathbb{R}^d$ to the population distribution \cite{mccann1995}. The convex function from this method does not have a formula and therefore cannot be directly analyzed as we have done for the K-function. However, we note that the center-outward quantile contours are proved to be \textit{connected}, closed, nested, with continuous boundaries \cite{hallin-etal-2021}. We foresee that the connectedness of the quantile contour would cause difficulties in detecting outliers for the non-standard distributions from \cite{xie2020local} depicted in Figure \ref{fig:Inverse-K-transform(Spiral)} and \ref{fig:Inverse-K-transform(SquareShape)} respectively, but we do not investigate this method any further in this article. 
          
%%%%%%%%%%%%%%%%%%%%%%%%%%%%%%%%%%%%%%%%%%%%%%%%%%
%%% K-transform and inverse K-transform
%%%%%%%%%%%%%%%%%%%%%%%%%%%%%%%%%%%%%%%%%%%%%%%
\section{The K-transform: Definition and Outlier Detection}\label{sec:K-transform}
Let $\left(\mathbb{R}^d,\mathcal{F}, \mu \right)$ be a probability
space. For any $s\in \mathbb{R}^d$, the Koltchinskii function (K-function) is
defined as the integral transform of $f(s,x)=\|s-x\|_2-\|x\|_2$ with respect to $\mu$:
\begin{equation}
	f_{\mu}\left(s\right) = \E f(s,x)=\int_{\mathbb{R}^d} \!\left(\|x-s\|_2-\|x\|_2\right)\,\mu(dx).\label{eq: Integral transform}
\end{equation}
It is easily verified that $f_{\mu}$ is a convex function. The K-transform (\cite{koltchinskii1997m}), is defined in as the sub-differential of $f_\mu$. For $s \in \mathbb{R}^d$,
\begin{equation}
	\partial f_{\mu}\left(s\right) = \int_{\mathbb{R}^d} \frac{s-x}{\|s-x\|_2} \,\mu(dx). \label{eqn: Koltchinskii transform}
\end{equation}
The K-transform is in one-to-one correspondence between $\mathbb{R}^d$
and the open unit ball $B_1(0)$. Now consider the conjugate function of $f_{\mu}$, given by, 
\begin{equation}\label{eq:M-quantile}
	f^{\ast}_{\mu}\left(v\right) = \underset{s \in \mathbb{R}^d}{\sup} \left[\langle s,v\rangle-f_{\mu}(s)\right], \quad  v\in B_1(0).
\end{equation}
The subdifferential of $f^{\ast}_{\mu}$ exists
(\cite{rockafellar1970convex}), and is given by $\partial
f^{\ast}_{\mu} = \left(\partial f_{\mu}\right)^{-1}$. The map
$\partial f^{\ast}_{\mu}$ is the inverse K-transform, and is referred
to as an M-quantile (\cite{koltchinskii1997m}). 
For dimension $d=1$, the K-transform is the signed cdf $\partial
f_{\mu}(s) = \mu\left(\left(-\infty,s\right]\right) -
\mu([s,\infty):=\skewF(s)$. In this case, the inverse K-transform is
interpreted as a quantile function. Also for $d=1$, there are
well-established techniques for outlier detection in data based on quantile plots, due to the facts that there is a natural ordering on the real line, and the derivative of the K-transform is the density. 

\begin{figure}[H]
	\centering \subfigure[]{ \label{fig:Inverse K-outliers}
		\includegraphics[width=0.40\textwidth]{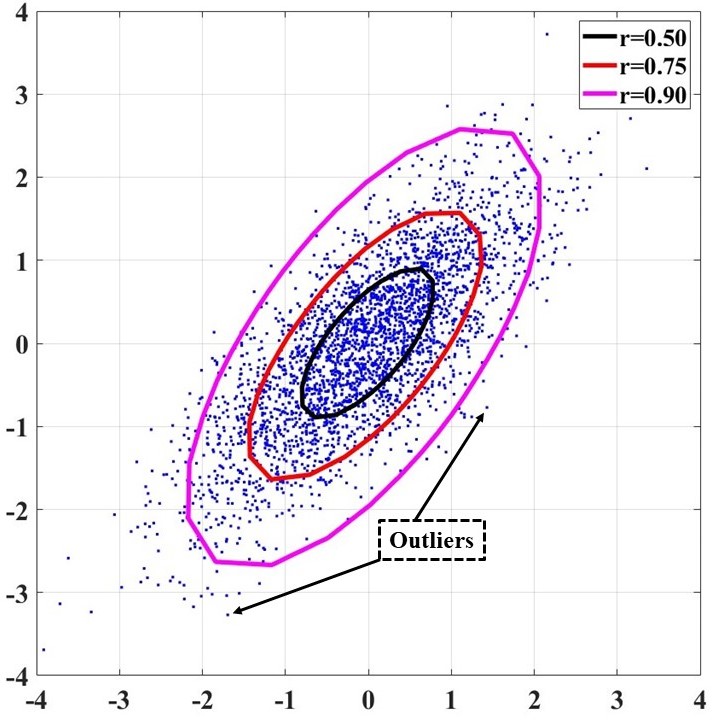}		
	} \subfigure[]{ \label{fig:K-transform}
		\includegraphics[width=0.41\textwidth]{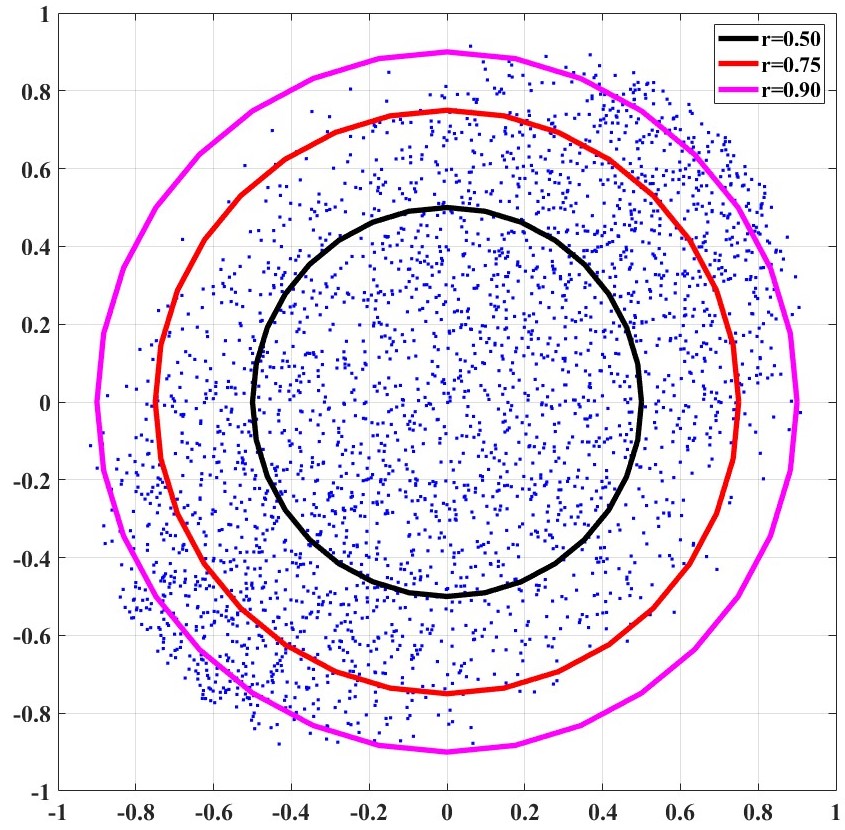}
	}
	\caption{K-transform and its inverse for a bivariate Gaussian with zero mean, unit standard
deviation, and correlation coefficient of $0.75$. The K-transform is
shown in Figure \ref{fig:K-transform}, and the inverse K-transform in
Figure \ref{fig:Inverse K-outliers}, along with 3000 data points simulated from the distribution.}
	\vspace{-0.05in}
	\label{fig:outliers}
\end{figure}

The inverse K-transform has been proposed as a means to detect outliers in
two-dimensional data (see e.g., \cite{chakraborty2001affine} and
\cite{chaouch2010design}). However, the underlying distributions were
chosen to be the Gaussian, which is one of the simplest of
cases. Figure \ref{fig:outliers} illustrates the potentially misleading result obtained in
\cite{chakraborty2001affine}. In Figure~\ref{fig:Inverse K-outliers}, we have plotted 3000 data points from a
bivariate normal distribution with zero mean, unit standard
deviation, and a correlation coefficient of $0.75$. The contours
here  are the inverse K-transform of the contours in Figure
\ref{fig:K-transform}. \cite{chakraborty2001affine} proposes  that
data points falling outside the inverse K-transform pertaining to
contours of the K-transform with a large radius (say
$r\geq 0.9$) may be classified as outliers.  
In Section \ref{sec:zoom-effect} we present numerical simulations
with more complex densities that reveals a very interesting
feature of the K-transform that we call the \textit{zoom-in}
effect. Due to this effect, we cannot interpret the radial distance
from the origin in Figure~\ref{fig:Inverse K-outliers} as representing
“quantiles”, a view that could lead to potentially erroneous outlier detection results. 

%%%%%%%%%%%%%%%%%%%
%%%%%%%%%%%%%%%%%%%

Our main result in the next section will show that for odd $d\geq 3$, a poly-Laplacian of the
K-function is proportional to the pdf. The proof of this result
requires the following non-trivial measure theoretic version of the
Leibniz Rule. We state this here as a major result in its own right.
\begin{theorem}[Lebesgue Measure Leibniz Rule] \label{th:Leibniz rule}
		Suppose $f(s,x)$ is bounded for all $x\in [a,b]$,
                and is both bounded and integrable with respect to
                Lebesgue measure $\lambda$ for all $s\in
                [\alpha,\beta]$.  For any Lipschitz continuous functions $g,h:[\alpha,\beta]
                \rightarrow [a,b]$, 
		let: \[ \phi(s) = \int_{h(s)}^{g(s)}\!f(s,x)\,dx. \]
		For all $s \in (\alpha,\beta)$ and for all $x \in
                [a,b] \setminus E_0,$ where $\lambda(E_0) = 0$,  let
                $\partial f(s,x)/\partial s$ be defined on
                $(\alpha,\beta) \times [a,b]$ and bounded in absolute
                value. Then, for almost all $s\in (\alpha,\beta)$,
                the derivative $d\phi(s)/ds$ exists and can be
                split as follows:
		\begin{equation} \frac{d}{ds}\phi(s) = f(s,g(s))\,\frac{dg(s)}{ds}  - f(s,h(s))\,\frac{dh(s)}{ds} + \int_{[h(s),g(s)] \setminus E_0}\!\frac{\partial}{\partial s}f(s,x)\;dx. \label{th:lebesgue} \end{equation}
\end{theorem}

%%%%%%%%%%%%%%%%%%%%%%%%%%%%%%%%%%%%%%%%%%%%%%%%%%%%%%%%%%%%%%%%%%%%%%%%%%%%%%
%%%%% OUR Results
%%%%%%%%%%%%%%%%%%%%%%%%%%%%%%%%%%%%%%%%%%%%%%%%%%%%%%%%%%%%%%%%%%%%%%%%%%%%%%%%%%%%%%%
\section{The K-transform in Relation to the Distribution}\label{sec:results}
In this section, we will present our main result concerning the
relationship of the K-transform to the distribution of the underlying
random vector $X$ in a probability
space $\left(\mathbb{R}^d,\mathcal{F}, \mu \right)$.  We will show that
for any odd values of the dimension $d$, it is in fact the poly-Laplacian of the
K-function that corresponds to
a scalar multiple of the pdf $\rho(x)$ of $X$.

Let $s \in \mathbb{R}^d$ be a given point. Let $\sigma_k$ denote the
hypersurface area measure on the $k$-sphere. For example, when $d=2$
we have
\[ d\sigma_1(\varphi_1)=d\varphi_1,d\sigma_2(\varphi_1,\varphi_2)=\sin(\varphi_2)\,d\varphi_1\,d\varphi_2,\]
and in general,
\[
  d\sigma_k(\varphi_1,\varphi_2,\cdots,\varphi_k)=\prod_{j=1}^k\,\sin^{j-1}(\varphi_j)\,d\varphi_j
,  \qquad k=1,2,3, \ldots.
\]
For simplicity, write $\Phi_k$ for the angular tuple
$(\varphi_1,\varphi_2, \cdots,\varphi_k)$. The volume element of
$x\in\mathbb{R}^d$ represented by the polar coordinates $\left(r,
  \varphi_1,\varphi_2,\cdots,\varphi_{d-1}\right)$, where $r$ is the
distance of $x$ from  $s$, is given by
\begin{equation*}
	dV(x) = r^{d-1}\,dr\,d\sigma_{d-1}(\Phi_{d-1}).
      \end{equation*}
An essentially bounded function $\rho: \mathbb{R}^d \rightarrow
\mathbb{R}_+$, satisfying
\[
  \displaystyle \int_{\mathbb{R}^d}\!\rho(x) \,dV(x) = 1,
\]
is said to be a \emph{bounded density} (bounded pdf). On the sigma-algebra of Borel subsets
of $\mathbb{R}^d$, we may use $\rho$ to define a Radon measure as
arising from a positive linear form on
$C_c\left({\mathbb{R}^d}\right)$, the set of continuous functions with
compact support on $\mathbb{R}^d$. Alternatively, one starts with an
absolutely continuous Radon measure with respect to the Lebesgue
measure, and takes its Radon-Nikodym derivative to be the  density
$\rho$.

\begin{remark}
Note that if $\mu(dx)=\rho(x)dV(x)$, then the K-transform is simply the
gradient of the K-function, so that $\partial f_{\mu}(s)=\nabla
f_{\mu}(s)$ in  \eqref{eqn: Koltchinskii transform}.
\end{remark}

%%%%%%%%% Lemma 01
The lemma below is a building block on which the theory is based. It
shows that the expectation of certain singular kernels is in fact well
defined. The only mild requirement is that the probability density
should be bounded. 
\begin{lemma}\label{Lemma 1}
	Suppose $\rho(x)$ is a bounded pdf for the random
        vector $X$ on $\mathbb{R}^d$, with $d\geq 3$. Then, for any
        $s\in\R^d$, the function defined as
\begin{equation} 
G(s) = \int_{\mathbb{R}^d
}\!\|s-x\|_2^{-k}\,\rho(x)\,dV(x), \label{eq: F(s)} 
\end{equation} 
exists for any $k = 0, \ldots, d-1.$
\end{lemma}

%%%%%%%%%%%%%%%%%%%%%%%%%%%%%%%%%%%%%%
%%%%%%% Lemma 02 
%%%%%%%%%%%%%%%%%%%%%%%%%%%%%%%%%%%%%%
The second lemma shows that the gradient of $G(s)$ defined above
exists for $0 \leq k \leq d-2.$ 

\begin{lemma}\label{Lemma 2}
Suppose $\rho(x)$ is a bounded pdf on $\mathbb{R}^d$, for $d \geq
3$. Then, for any $s \in \mathbb{R}^d$ and $k = 0,\ldots,d-2$,
	the Frechet derivative of $G(s)$ defined in \eqref{eq: F(s)} with respect to $s$, exists and is given by 
	\begin{equation} \nabla_{s}\,G(s) = -k\,\int_{\mathbb{R}^d }\!\frac{(s-x)}{\|s-x\|_2^{k+2}}\,\rho(x)\,dV(x).\end{equation}
\end{lemma}

%%%%%%%%%%%%%%%%%%%%%%%%%%%
%%%%%% Lemma 03
%%%%%%%%%%%%%%%%%%%%%%%%
\vspace{0.2in}
The third lemma shows that the Hessian of $G(s)$ exists for $0 \leq k \leq d-3.$ The Hessian computation is important as the Laplacian may be computed as the trace of the Hessian matrix.

\begin{lemma}\label{Lemma 3}
	Suppose $\rho(x)$ is a bounded pdf on $\mathbb{R}^d$ for
        $d \geq 3$. Then, for any $s \in \mathbb{R}^d$ and $k=0,\ldots,d-3$, 
	the Hessian of $G(s)$ defined in \eqref{eq: F(s)} with respect to $s$, is given by: 
	\begin{equation}
		D_sG(s) = -k \int_{\mathbb{R}^d} \left(\frac{I_{d}}{\|s-x\|_2^{k+2}}-\left(k+2\right) \frac{\left(s-x\right)\left(s-x\right)^{t}}{\|s-x\|_2^{k+4}}\right)\rho\left(x\right)dV\left(x\right). \label{eqn: Jacobian of F(s)}
	\end{equation}  
\end{lemma}

We are now ready to produce results for the K-transform. The lemma
below considers the univariate case ($d =1$); the special case of odd
dimension that is omitted in the statement of Theorem
\ref{Theorem}. The result is proved in
\cite{koltchinskii1997m} using convex analysis. Here we
use real analysis techniques, because for dimensions greater than one, they
allow one to obtain expressions for the poly-Laplacian of the K-function. 

%%%%%%%%%%% Lemma 04 %%%%%%%%%%%%%%%
\begin{lemma} \label{Lemma 4}
	Consider the probability space $\left(\mathbb{R}, \mathcal{F},
          \mu\right)$, where the measure $\mu$ of the
        random element $X$ has a bounded density, $\rho(x)$, on
        $\mathbb{R}$. Consider the $d=1$ version of the K-function
        defined in \eqref{eq: Integral transform}:
        \[
f_\mu(s) = \int_{\R}(|s-x|-|x|)\rho(x)dx.
        \]
        Then, for almost every $s \in \mathbb{R}$, the second derivative of
$f_{\mu}(s)$ exists, and is given by, 
	\begin{equation}
		f_{\mu}''(s)= 2\rho(s). \label{eqn: Lemma 4.2}
	\end{equation}
\end{lemma}

%%%%%%%%%%%%%%%%%%%%%%%%%%%%%%%%%%%%%%%%%%%
%%%%%% Theorem Starts 
%%%%%%%%%%%%%%%%%%%%%%%%%%%%%%%%%%%%%%%
With the above lemmas in place, we can now state the main theoretical result of this article. Taken together, we show that for
odd dimensions $d \geq 3$, the poly-Laplacian of order $(d+1)/2$ of the
K-function is equal to a scalar multiple of the density (almost
everywhere), with the scalar depending only on $d$. As regions of low
probability density are commonly associated with outliers and
anomalies in data (see
\cite{chandola2009,aggarwal2015}), this result shows
that, contrary to emerging belief, it is not the K-transform $\nabla
f_{\mu}(s)$ that is directly related to the density, but the
poly-Laplacian $\Delta^{(d+1)/2}f_{\mu}(s)$ instead. Hence, it is
dubious to use the K-transform to compute “quantiles” in odd
dimensions $d \geq 3$; a situation that is completely different from
the $d =  1$ case.

Using Lemma \ref{Lemma 3}, we can obtain a
different type of result for even dimensions. However, we can only
establish a lower bound for $\|\nabla \, f_\mu(s)\|_2$  for $d = 2$, and
an upper bound for $\|\nabla \,\Delta^{(d-2)/2} f_\mu(s)\|_2$ for
$d \geq 4$, in terms of $\rho(s)$. Unfortunately  these results do not
appear to have as clear an interpretation as those for odd dimensions.

\begin{theorem} \label{Theorem}
	Consider the probability space $\left(\mathbb{R}^d,
          \mathcal{F}, \mu\right)$, where the probability measure
        $\mu$ has a bounded density $\rho(x)$ in $\mathbb{R}^d$, and whose
        K-function $f_{\mu}(s)$ is as defined in \eqref{eq: Integral transform}. 
Let $\Delta^{j}f_{\mu}(s)$ be the poly-Laplacian of order $j$ of $f_{\mu}(s)$,
defined as:
\begin{equation}\label{eqn: Theorem}
\Delta^{j}f_{\mu}(s) = (-1)^{j+1} (d-1)(d-3)\dots (d-2j+1) \, \E\left(\frac{1}{\|s-X\|_2^{2j-1}}\right).
\end{equation}
Then, for any $s \in \mathbb{R}^d$ and odd values of the dimension
        $d\geq 3$, $\Delta^{j}f_{\mu}(s)$ exists for all
        $j=1,\ldots,N$, where $N=(d-1)/2$.
\end{theorem}
\begin{corollary} \label{corollary}
	For any odd $d\geq 3$, and for almost every $s\in
        \mathbb{R}^d$, we have that: 
	%%%%%%%%%%%%%%%%%%%%%%%%%%%%%
	\begin{equation}\label{eqn: Corollary 1.1}
		\Delta^{N+1} f_{\mu}(s) = \alpha\, \rho(s), 
	\end{equation}
	%%%%%%%%%%%%%%%%%%%%%%%%%%%%%
	where $ \alpha = (-1)^{N-1}\,(d-1)\,\cdots\,(2) \, d\,(d-2)\pi^{d/2}/\Gamma\left(\frac{d}{2}+1\right) $, and $\Gamma$ is the gamma function.
\end{corollary}
%%%%%%%%%%%%%%%%%%%%%%%%%%%%%%%%%%%%%%%%%%%%%%%%%%%%%%%%%%%%%%%%
%%% Numerical results for irregular densities on $\mathbb{R}^d
%%%%%%%%%%%%%%%%%%%%%%%%%%%%%%%%%%%%%%%%%%%%%%%%%%%%%%%%%%%%%%%%%

Note that the first few poly-Laplacians are (assuming odd $d > 5$):
\[
\Delta_{s} f_{\mu}(s) = (d-1) \, \E\left(\frac{1}{\|s-X\|_2}\right), \qquad
		\Delta_{s}^2f_{\mu}(s) = (-1)\,(d-1)\,(d-3) \,
                \E\left(\frac{1}{\|s-X\|_2^3} \right),
              \]
and
              \[
\Delta_{s}^3f_{\mu}(s) = (-1)^2\,(d-1)\,(d-3)\,(d-5)\, \E\left(\frac{1}{\|s-X\|_2^5} \right) .
\]
The final polyharmonic is:
\begin{equation}
  \Delta_{s}^{N}f_{\mu}(s) = (-1)^{N-1}\,(d-1)\,\cdots\,(2)\,
  \E\left(\frac{1}{\|s-X\|_2^{d-2}} \right)  . \label{eq: final
    polyharmonic}
\end{equation}

%%%%%%%%%%%%%%%%%%%%%%%%%%%%%%%%%%%%%%%%%%%%%%%%%%%%%%%%%%%%%%%%% 
\section{The zoom-in effect for $\mathbb{R}^2$} \label{sec:zoom-effect}
%%%%%%%%%%%%%%%%%%%%%%%%%%%%%%%%%%%%%%%%%%%%%%%%%%%%%%%%%%%%%%%%%
The inverse K-transform is useful for detecting outliers in simple
distributions on $\mathbb{R}^2$, for instance the Gaussian, as demonstrated by \cite{chakraborty2001affine}.  
However, \cite{chen2008} cautioned  that this claim is incorrect for
non-standard distributions in $\mathbb{R}^2$ such as a mixture of
Gaussians. In this section, we demonstrate why this happens by
highlighting a special feature of the K-transform that we call the
“zoom-in” effect. 

\begin{figure}[H]
	\centering \subfigure[]{ \label{fig:Inverse-K-Transform (BananaShape)}
		\includegraphics[width=0.40\textwidth]{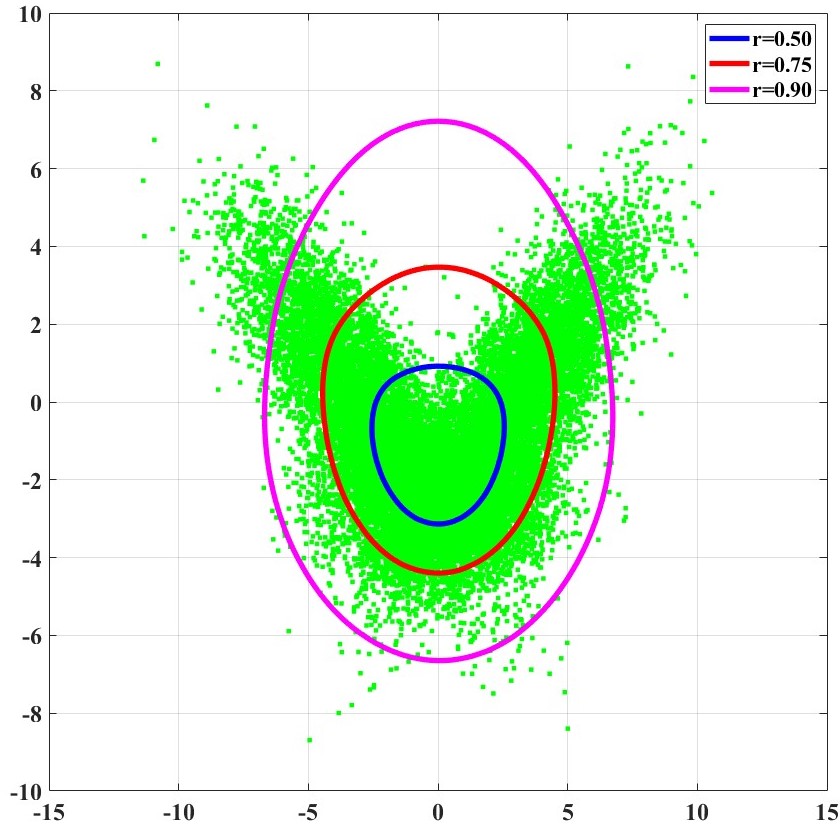}		
	} \subfigure[]{ \label{fig:K-Transform (BananaShape)}
		\includegraphics[width=0.40\textwidth]{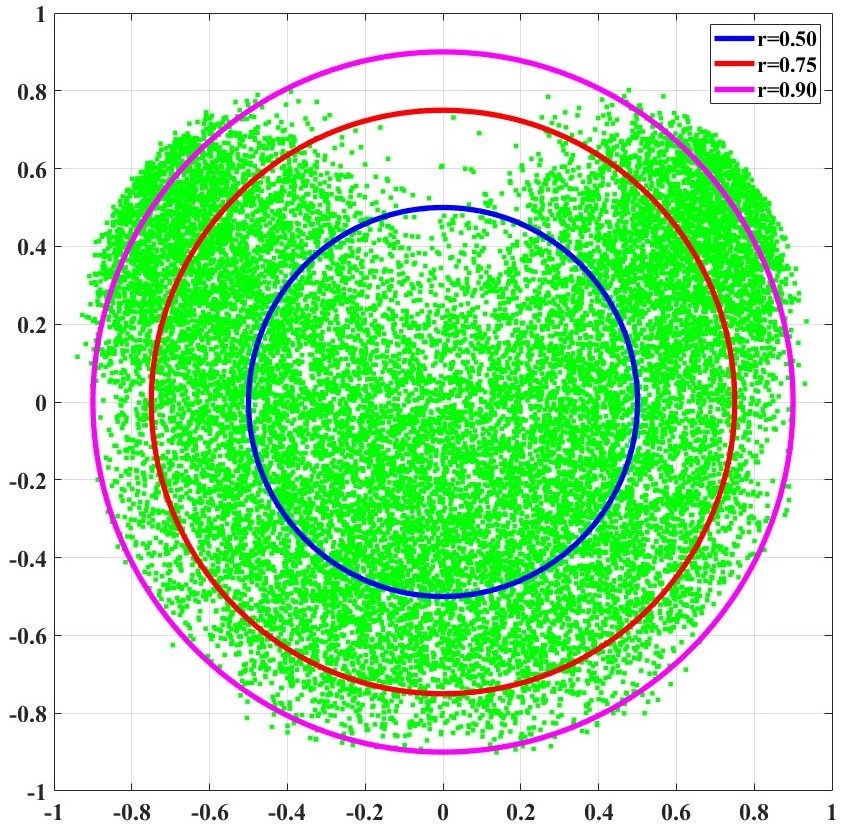}
	}
	\caption{Banana-shaped distribution with $n=20,000$ sample points \cite{hallin-etal-2021}. Concentric circles of radius 0.50 (blue), 0.75 (red), and 0.90 (maroon) are shown in Figure \ref{fig:Inverse-K-Transform (BananaShape)}. Their corresponding K-transforms are in Figure \ref{fig:K-Transform (BananaShape)}. This example illustrates the zoom-in effect of the K-transform. Points outside the support of the distribution in the domain occupy a small region close to the boundary of the unit circle in the co-domain. Hence, $1-\|\nabla f\|_2$ cannot be interpreted as a depth function. }
	\vspace{-0.05in}
	\label{fig:BananaShape}
\end{figure}

\begin{figure}[H]
	\centering \subfigure[]{ \label{fig:Inverse-K-transform(Spiral)}
		\includegraphics[width=0.40\textwidth]{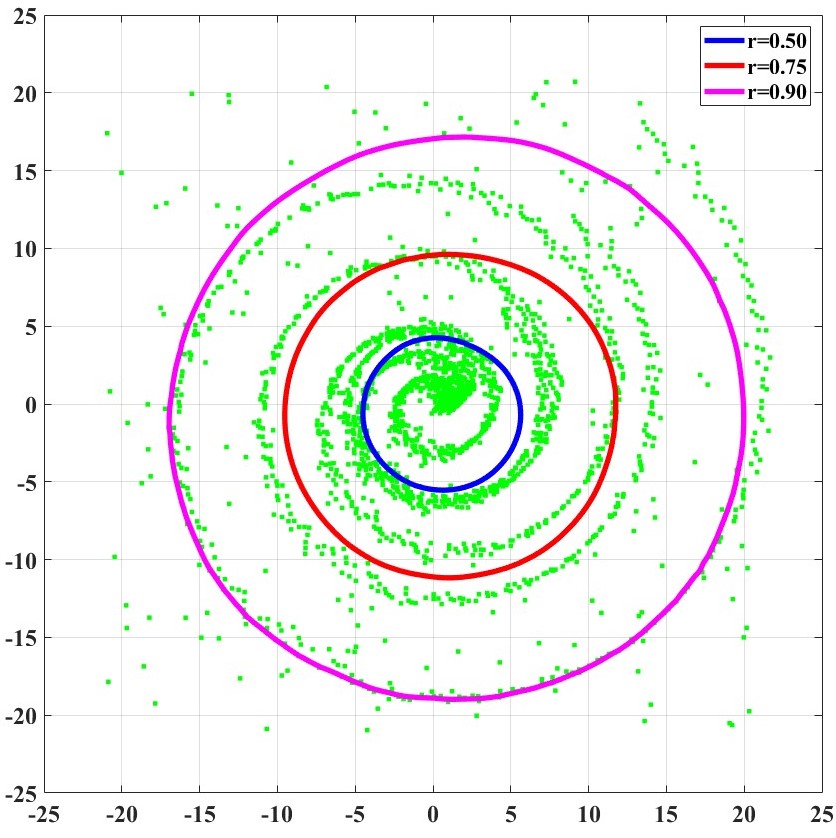}		
	} \subfigure[]{ \label{fig:K-transform(Spiral)}
		\includegraphics[width=0.41\textwidth]{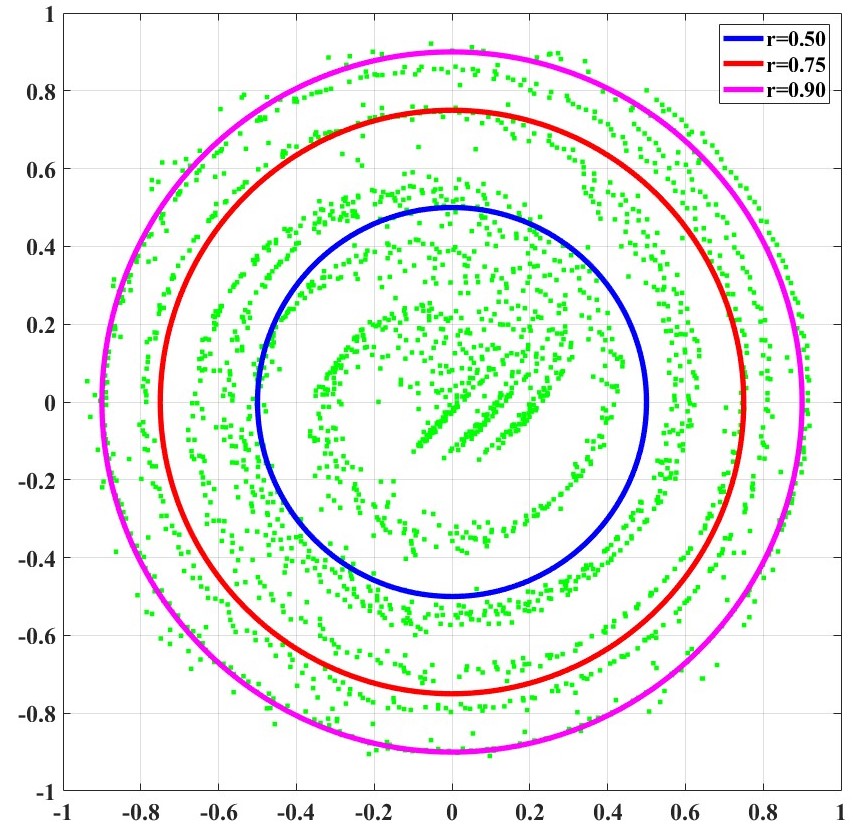}
	}
	\caption{Spiral-shaped distribution with $n=2259$ sample points \cite{xie2020local}. Figure \ref{fig:K-transform(Spiral)} shows the K-transform and \ref{fig:Inverse-K-transform(Spiral)} shows the inverse K-transform. This example illustrates zoom-in effect and the low power of the K-transform method. Points in the empty spaces are within the $r=0.5$ M-quantile while points in the outer spiral of the data are outside the circle of radius $0.9$ in the co-domain.}
	\vspace{-0.05in}
	\label{fig:Spiral}
\end{figure}

\begin{figure}[H]
	\centering \subfigure[]{ \label{fig:Inverse-K-transform(SquareShape)}
		\includegraphics[width=0.40\textwidth]{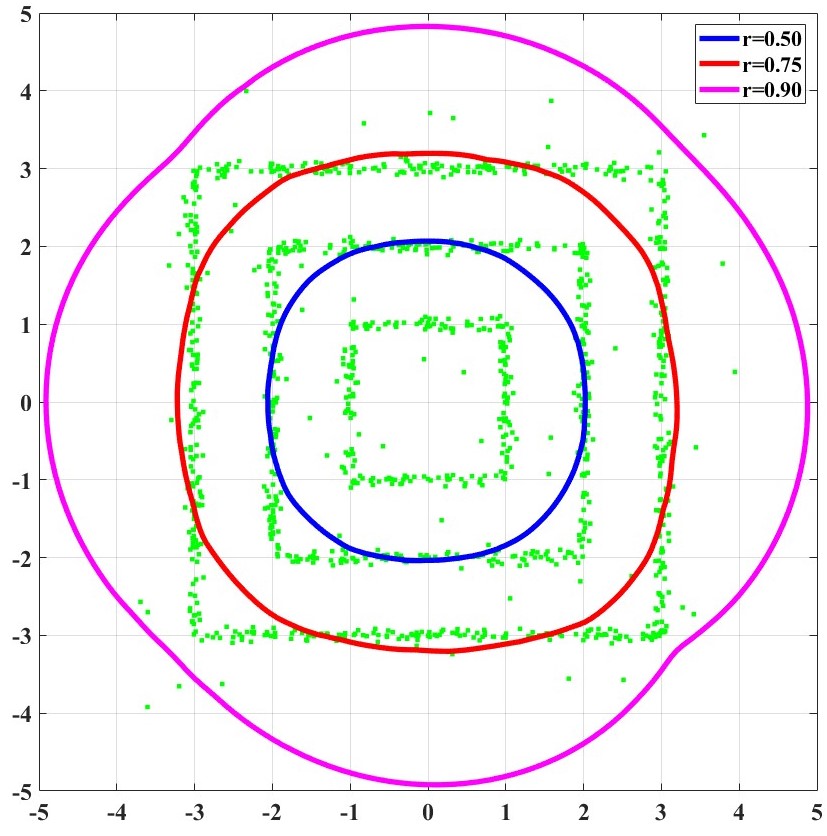}		
	} \subfigure[]{ \label{fig:K-transform(SquareShape)}
		\includegraphics[width=0.41\textwidth]{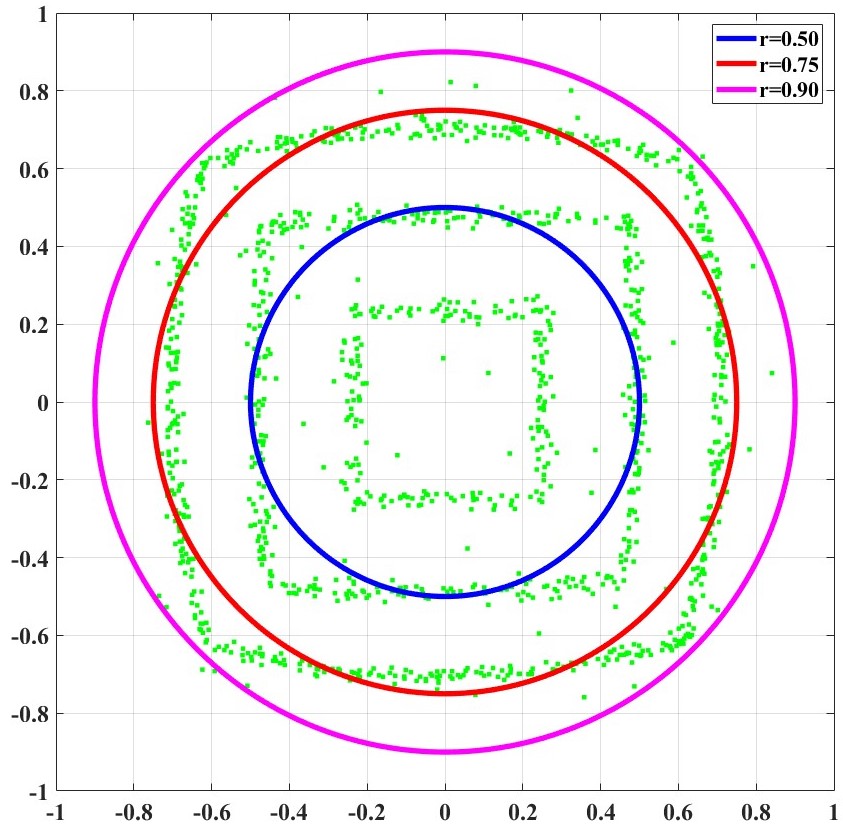}
	}
	\caption{This example illustrates zoom-in effect and the low power of the K-transform method for a non-standard square-shaped distribution with $n=1242$ sample points \cite{xie2020local}. Figure \ref{fig:K-transform(SquareShape)} shows the K-transform and \ref{fig:Inverse-K-transform(SquareShape)} shows the inverse K-transform. Points in the empty spaces within the $r=0.5$ M-quantile should be outliers but are not classified as such by the method. }
	\vspace{-0.05in}
	\label{fig:SquareShape}
\end{figure}

Consider the unit disk in $\mathbb{R}^2$ centered at the origin with a
uniform density of $\pi^{-1}$.  By symmetry of the disk, the
K-transform at the point $s=(s_1,s_2)$ is given by: 
%%%%%%% Equation  %%%%%%%%%%%%%%%%%%%%%%
\begin{equation}
	\nabla f_{\mu}(s_1,s_2) = \int_{B_1(0)}\! \frac{(s_1,s_2)-(x_1,x_2)}{\|(s_1,s_2)-(x_1,x_2)\|_2} \,\mu(dx).   \label{eq:main equation} 
\end{equation}
%%%%%%%%%%%%%%%%%%%%%%%%%%%%%%%%%%%%%%%%
As discusssed in Section~\ref{sec:K-transform}, the K-transform will
map all of $\mathbb{R}^2$ onto the unit disk in the codomain. Figures \ref{fig:BananaShape}, \ref{fig:Spiral}, and \ref{fig:SquareShape} depict the zoom-in effect on the banana-shaped \cite{hallin-etal-2021}, spiral-shaped \cite{ha2015precise,xie2020local}, and square-shaped \cite{ha2015precise,xie2020local} distributions. Note that the \emph{zoom-in} effect refers to the tendency of the K-transform to magnify regions of the support, so that they occupy a disproportionately larger area of the unit disk.

 The banana-shaped distribution is a mixture of three Gaussian distributions with $n=20,000$ sample points, given by,  

\begin{equation*}
 \frac{3}{8} N(-3\tau_{1},\Sigma_1)+\frac{3}{8}N(3\tau_1,\Sigma_2)+\frac{1}{4} N(-\frac{5}{2}\tau_2,\Sigma_3),   
\end{equation*}
where $\tau_1 = \begin{pmatrix}
	1 \\
	0 
\end{pmatrix}$, 
$\tau_2 = \begin{pmatrix}
	0 \\
	1 
\end{pmatrix}$,
$\Sigma_1= \begin{pmatrix}
	5 & -4 \\
    -4 & 5 
\end{pmatrix}$,
$\Sigma_2= \begin{pmatrix}
	5 & 4 \\
	4 & 5 
\end{pmatrix}$, and, 
$\Sigma_3= \begin{pmatrix}
	4 & 0\\
	0 & 1 
\end{pmatrix}$. 
The contour lines in Figure \ref{fig:Inverse-K-Transform (BananaShape)}, \ref{fig:Inverse-K-transform(Spiral)}, and \ref{fig:Inverse-K-transform(SquareShape)} are the inverse K-transform of the contours with radial distances of 0.5 (blue), 0.75 (red), and 0.90 (maroon) in the codomain. We have used our developed MM-type algorithm to compute the inverse K-transform of the contours that are in the codomain. \cite{chakraborty2001affine, girard2017intriguing} described the inverse K-contour lines as quantile contours and \cite{chakraborty2001affine} classified points that are outside of radial distance 0.9 as outliers. From these figures, we can observe that because of the zoom-in effect, inlier points become close to the boundary of the contour lines with a radial distance of 0.9, and outlier points get close to the inliers.

%effect refers to the tendency of the K-transform to magnify regions of
%the support, so that they occupy a disproportionately larger area of
%the unit disk. Chakraborty called the inverse K-contour lines quantile contours and classified points that are outside of radial distance 0.9 as outliers. From these figures, we can observe that because of the zoom-in effect, inlier points become close to the boundary of the contour lines with a radial distance of 0.9, and outlier points get close to the inliers.
%Figures
%\ref{fig:Cross} and \ref{fig:jellyBean} depict the outcome of this
%transformation in two particular instances, where uniform measures
%are distributed over a ``cross'' and a ``bean'' shape. The \emph{zoom-in}
%effect refers to the tendency of the K-transform to magnify regions of
%the support, so that they occupy a disproportionately larger area of
%the unit disk. This is particularly evident in Figure~\ref{fig:jellyBean}.

%For instance, in Figure \ref{fig:Cross}, we can see that the
%%%distance of $A'$ from $B'$ is larger than the distance of $C'$ from
%$B'$. Equating distance from the origin with increasing likelihood of
%outlyingness, would lead to the conclusion that $A$ is more likely to be
%an outlier than $C$. However, $A$ is actually on the boundary of the
%support, while $C$ is well outside! A similar effect is observed in Figure~\ref{fig:jellyBean}.

%%%%%%%%%%%%%%%%%%%%%%%%%%%%%%%%%%%%%%
%% Banana-Shaped
%%%%%%%%%%%%%%%%%%%%%%%%%%%%%%%%%%%%%%
%%%%%%%%%%%%%%%%%%%%%%%%%%%%%%%%%%%%%%%%%%%
%% Spiral-Shaped 
%%%%%%%%%%%%%%%%%%%%%%%%%%%%%%%%%%%%%%%%%%%%%

%%%%%%%%%%%%%%%%%%%%%%%%%%%%%%%%%%%%%%%%%%%
%% Squared-Shaped 
%%%%%%%%%%%%%%%%%%%%%%%%%%%%%%%%%%%%%%%%%%%%%

%%%%%%%%%%%%%%%%%%%%%%%%%%%%%%%%%%%%%%%%%%%%%%%%%%%%%%%%%%%%%
%%%%%% Conclusion
%%%%%%%%%%%%%%%%%%%%%%%%%%%%%%%%%%%%%%%%%%%%%%%%%%%%%%%%%%%%%%
\section{Conclusion}
Regions of low probability density are commonly associated with
outliers and anomalies in data. In this article, we have  shown that
for multivariate data on spaces with odd dimension $d$ greater
than or equal to three, it is not the K-transform $\nabla f_\mu(s)$
that is directly proportional to the probability density, but the
poly-Laplacian $\Delta^{(d+1)/2}\,f_\mu(s)$ instead. This leads to the
natural conclusion that the K-transform is not a useful object from
which to attempt to compute “quantiles” in higher dimensions; a
situation that is very different from the one-dimensional case. For
the lowest even dimensional case of $d = 2$, we elucidate a zoom-in
feature of the K-transform whereby regions of higher probability
density are magnified when mapped to the co-domain. Due to this
effect, it may not be appropriate to think of contours of equal radial
distance from the origin in the co-domain as corresponding to
quantiles. In summary, the established lore of interpreting the
K-transform as a distribution and its inverse as a quantile function, applies
only to the one-dimensional case.

\begin{appendix}
	\section{Algorithm for Computation of Geometric Quantiles}\label{sec:app-geomQ-algo}
	In the following steps, we describe an MM type
        algorithm (\cite{lange2016mm}) for the computation of the geometric quantile in \eqref{eqn:geom-quantile}.
\begin{description}
\item[Step 0:] Randomly choose an initial value of $s$, say $s_0$. \\
\item[Step 1:] Specify the surrogate objective function to be 
	%%%%%%%%%% Equation 10 %%%%%%%%%%%%%
	\begin{equation}\label{eqn: MM 1.2}
		\underset{s \in \mathbb{R}^d}{\arg \min} \;\frac{1}{2} \E \left(\frac{\|s-X\|_{2}^2}{\|s_0-X\|_{2}}+\|s_0-X\|_{2} - 2\|X\|_{2}\right)-\langle s,v \rangle, v \in B_1(0).  
	\end{equation}
	%%%%%%%%%%%%%%%%%%%%%%%%%%%%%%%%%%%
\item[Step 2:] Compute the gradient of (\ref{eqn: MM 1.2}) as, 
	\begin{equation}
		v=\E \left(\frac{s-X}{\|s_0-X\|_2}\right).
	\end{equation}
\item[Step 3:] For  $k\geq 1$, iteratively update the value of $s$ until convergence is reached.
	%%%%%%%%%% Equation 11 %%%%%%%%%%%%%
	\begin{equation}
		s_{k}  =\frac{v+\E\left(X\|s_{k-1}-X\|_2^{-1}\right)}{\E\left(\|s_{k-1}-X\|_2^{-1} \right)}.
              \end{equation}
\end{description}

%%%%%%%%%%%%%%%%%%%%%%%%%%%%%%%%%%%% 
	%%%% Proof of Leibnitz Rule
%%%%%%%%%%%%%%%%%%%%%%%%%%%%%%%%%%%%
\section{Proof of Theorem \ref{th:Leibniz rule}}
\begin{proof}%[Proof of Theorem \ref{th:Leibniz rule}]
Suppose $dg(t)/dt$ and $dh(t)/dt$ exist for all $t \in
(\alpha,\beta),$ and let $\{t+t_k\}_{k=1}^\infty$ be any sequence
contained in $(\alpha,\beta)$ converging to $t.$ We then have:
\begin{align}  \frac{\phi(t + t_k)- \phi(t)}{t_k} & =  \frac{1}{t_k}\left[\int_{h(t+t_k)}^{g(t+t_k)}\!\!\!f(x,t+t_k)\,dx - \int_{h(t)}^{g(t)}\!\!\!f(x,t)\,dx \right]  \nonumber \\
	& =  \frac{1}{t_k}\left[\int_{h(t+t_k)}^{g(t+t_k)}\!\!\!f(x,t+t_k)\,dx - \int_{h(t+t_k)}^{g(t)}\!\!\!f(x,t+t_k)\,dx \right] + 
	\nonumber\\ &  \frac{1}{t_k}\left[ \int_{h(t+t_k)}^{g(t)}\!\!\!f(x,t+t_k)\,dx - \int_{h(t)}^{g(t)}\!\!\!f(x,t+t_k)\,dx \right] + 
	\nonumber \\ &  \frac{1}{t_k}\left[\int_{h(t)}^{g(t)}\!\!\!f(x,t+t_k)\,dx- \int_{h(t)}^{g(t)}\!\!\!f(x,t)\,dx \right] \nonumber \\
	& =  \frac{1}{t_k}\left[ \int_{g(t)}^{g(t+t_k)}\!\!\!f(x,t+t_k)\,dx - \int_{h(t)}^{h(t+t_k)}\!\!\!f(x,t+t_k)\,dx \right] \nonumber \\ &  + \frac{1}{t_k}\left[\int_{h(t)}^{g(t)}\!\!\![f(x,t+t_k)\,dx- f(x,t)]\,dx \right]. \label{three terms}
\end{align} 
Let $|f(t,x)| \leq M$ on $[\alpha,\beta] \times [a,b]$ and
$\vert\partial f(t,x)/\partial t\vert\leq K $ on $(\alpha,\beta)
\times (a,b)$. We note that for $\varepsilon > 0$ , there exists
$\delta > 0$ such that for all $ |t_k| < \delta$ and $x \in[a,b]$, $ |
f(x,t+t_k) - f (x,t)| < K \,|t-t_k|$. In addition, since $g$ is
Lipschitz continuous, we have that
\begin{align*} \Big{\vert} \frac{1}{t_k}\,\int_{g(t)}^{g(t+t_k)}\!\!\!
  \left(f(x,t+t_k) - f(x,t)\right)\,dx \Big{\vert} & \leq \Big{\vert}
                                                     \frac{1}{t_k}\,\int_{g(t)}^{g(t+t_k)}\!\!\!|
                                                     f(x,t+t_k) -
                                                     f(x,t)|\,dx \Big{\vert}  \\ 
& \leq  \Big{\vert} \frac{K}{t_k}\,\int_{g(t)}^{g(t+t_k)}\!\!\!|t+t_k - t| \,dx \Big{\vert} \\ & \leq K\,|g(t+t_k) - g(t)| \\ & \leq K\, \|g\|_\infty \,|t_k|. \end{align*}
Therefore,\;\;\; 
\begin{equation} \lim_{t_k \rightarrow 0}\, \Big{\vert} \frac{1}{t_k}\,\int_{g(t)}^{g(t+t_k)}\!\!\! [f(x,t+t_k) - f (x,t)]\,dx \Big{\vert} = 0. \label{key limit} \end{equation}
Now, by the change of variable formula (see e.g., \cite[Sec.~38.3]{mcshane-botts}):
%(see e.g., \cite[Sec.~17]{billingsley1995}):
\begin{align} \frac{1}{t_k}\int_{g(t)}^{g(t+t_k)}\!\!\!\!f(x,t)\,dx 
	%	&=  \frac{1}{t_k}\int_t^{t+t_k}\!\!\!\!f(g(\tau),t+t_k)\,g'(\tau)\,d\tau \\
	&=
   \frac{1}{t_k}\int_t^{t+t_k}\!\!\!\!f(g(\tau),t)\,g'(\tau)\,d\tau.
\end{align}
As this integrand is integrable, we then have by the
Lebesgue Differentiation Theorem that for almost every $t \in [\alpha,\beta]$,
%\begin{equation}  \lim_k \,\frac{1}{t_k}\int_t^{t+t_k}\!\!!f(g(\tau),t+t_k)\,g'(\tau)\,d\tau = \lim_k \,\int_0^1\!\!\!F_k(\sigma)\,d\sigma = \int_0^1\!\!\!\lim_k \,F_k(\sigma)\,d\sigma = f(g(t),t)\,g'(t). \label{first term} \end{equation}
\begin{equation}  \lim_k \,\frac{1}{t_k}\int_t^{t+t_k}\!\!\!f(g(\tau),t)\,g'(\tau)\,d\tau  = f(g(t),t)\,g'(t). \label{approx first term} \end{equation}
By \eqref{key limit} and \eqref{approx first term},
\begin{align}
	\lim_k \,\frac{1}{t_k}\,\int_{g(t)}^{g(t+t_k)}\!\!\! f(x,t+t_k) \,dx	&= \lim_k \,\frac{1}{t_k}\,\int_{g(t)}^{g(t+t_k)}\!\!\! [f(x,t+t_k) - f(x,t)]\,dx  \nonumber \\
	& \quad - \lim_k \,\frac{1}{t_k}\,\int_{g(t)}^{g(t+t_k)}\!\!\! f(x,t)\,dx \nonumber \\
	&=   f(g(t),t)\,g'(t). \label{first term}
\end{align}
By the same approach, the second integral in (\ref{three terms}) can be shown to be for almost every $t \in [\alpha,\beta]$
\begin{equation} \lim_k \,\frac{1}{t_k}\int_t^{t+t_k}\!\!\!f(h(\tau),t+t_k)\,h'(\tau)\,d\tau = f(h(t),t)\,h'(t). \label{second term} \end{equation}
Consider now the final term in (\ref{three terms}). Let $x\notin
E_0$. By the Mean Value Theorem applied to the interval $[t,t+t_k]$ or
$[t+t_k,t]$, depending on the sign of $t_k$, we then have that:
%: \[ \int_{h(t)}^{g(t)}\!\!\!\frac{f(x,t+t_k)\,dx- f(x,t)}{t_k}\,dx.\] By the mean value theorem, for all $x \in [a,b] \setminus E_0$: 
\[
  \Big{\vert}\frac{f(x,t+t_k)\,dx- f(x,t)}{t_k}\Big{\vert} =
  \Big{\vert}\frac{\partial f}{\partial t}(\bar{t},x)\Big{\vert} \leq
  K,
\]
for some $\bar{t}$ in the segment joining $t$ and $t+t_k$. Now, invoking the
(Lebesgue) Dominated Convergence Theorem, we obtain:
\begin{eqnarray} \lim_k \,\int_{h(t)}^{g(t)}\!\!\!\frac{f(x,t+t_k)\,dx- f(x,t)}{t_k}\,dx &=& \int_{h(t)}^{g(t)}\!\!\!\lim_k \,\frac{f(x,t+t_k)\,dx- f(x,t)}{t_k}\,dx \nonumber\\ &=& \int_{[h(t),g(t)] \setminus E_0}\!\!\!\frac{\partial f}{\partial t}(\bar{t},x)\,dt. \label{third term}\end{eqnarray}
Collecting \eqref{first term} -- \eqref{third term}, we obtain finally  that for almost every $t \in (\alpha,\beta)$:
\begin{eqnarray} \nonumber \frac{d}{dt}\phi(t)& = & \lim_k\,\frac{1}{t_k}\left[\int_{h(t+t_k)}^{g(t+t_k)}\!\!\!f(x,t+t_k)\,dx - \int_{h(t)}^{g(t)}\!\!\!f(x,t)\,dx \right]  \\
	& = & f(t,g(t))\,\frac{dg(t)}{dt}  - f(t,h(t))\,\frac{dh(t)}{dt} + \int_{[h(t),g(t)] \setminus E_0}\!\frac{\partial}{\partial t}f(t,x)\;dx.  \end{eqnarray}
\end{proof}

%%%%%%%%%%%%%%%%%%%%%%%%%%%%%%%%%%%%
	%%%% Proof of Lemma 3.1
%%%%%%%%%%%%%%%%%%%%%%%%%%%%%%%%%%%%
	\section{Proof of Lemma \ref{Lemma 1}}
	
	%\begin{lemma}
	%	Suppose $\rho$ is a bounded density on $\mathbb{R}^d$ for $d \geq 3$. For $s \in \mathbb{R}^d$,  \begin{equation} \int_{\mathbb{R}^d}\!\|s-x\|_2^{-k}\,\rho(x)\,dV(x) < \infty,\end{equation} for $k = 0, \cdots, d-1.$
	%\end{lemma}
	\begin{proof}
		Set $z = x-s$. Denote $r = \|z\|$. We will use the volume form for $\mathbb{R}^d$ in the coordinates $z$.
		\begin{align*}
			\int_{\mathbb{R}^d}\!\|s-x\|_2^{-k}\,\rho(x)\,dV(x) &= \int_{B_1(0)}\!\|z\|_2^{-k}\,\rho(z+s)\,dV(z) + \int_{B_1^c(0)}\!\|z\|_2^{-k}\,\rho(z+s)\,dV(z)\\
			&= \int_{B_1(0)}\!\!\!\!\!\!\!\!\!\rho(z+s)\,r^{d-1-k}\,dr\,d\sigma_{d-1}(\Phi_{d-1}) +  \int_{\mathbb{R}^d\setminus B_{1}(0)}\!\!\!\!\!\!\!\!\!\!\!\!\!\!\!\!\rho(z+s)\,r^{d-1-k}\,dr\,d\sigma_{d-1}(\Phi_{d-1}).
		\end{align*}
		As $k \geq 0$, \[ r^{d-1-k} \leq \left\{\begin{array}{rl} 1; & r \leq 1 \\ r^{d-1}; & r > 1. \end{array} \right.  \] Hence,
		\begin{align*}
			\int_{\mathbb{R}^d}\!\|s-x\|_2^{-k}\,\rho(x)\,dV(x) &\leq \int_{B_1(0)}\!\!\!\!\rho(z+s)\,dr\,d\sigma_{d-1}(\Phi_{d-1}) +  \int_{\mathbb{R}^d\setminus B_{1}(0)}\!\!\!\!\!\!\!\!\!\!\!\!\!\rho(z+s)\,r^{d-1}\,dr\,d\sigma_{d-1}(\Phi_{d-1}) \\
			&\leq \|\rho\|_\infty\,\int_{B_1(0)}\!\!\!\,dr\,d\sigma_{d-1}(\Phi_{d-1}) +  \int_{\mathbb{R}^d}\!\!\rho(z+s)\,r^{d-1}\,dr\,d\sigma_{d-1}(\Phi_{d-1}) \\
			&<  C + 1 \\
			& < \infty.
		\end{align*}
	\end{proof}
	
	%%%%%%%%%%%%%%%%%%%%%%%%%%%%%%%%%%%%%%%%%%%%%%%%
	%%%% Proof of Lemma 3.2
	%%%%%%%%%%%%%%%%%%%%%%%%%%%%%%%%%%%%%%%%%%%%%%%%%
	\section{Proof of Lemma \ref{Lemma 2}}
	\begin{proof}
Write $\|\cdot\|$ to mean the usual
                Euclidean norm $\|\cdot\|_2$ in order to simplify the notation in
                the rest of the proof.          There is a need for subtlety because the integrands of $G$ and its proposed gradient are undefined for $x= s$. We need to excise a small region around the point $s$ in the variable $x$  while dealing with the integrals.
		To this end, let $R > 0$. Since the integrand of $G(s)$ given by $g(s,x) =
                \|s-x\|^{-k}$ is a differentiable and bounded function
                on $\mathbb{R}^d \setminus B_R(s)$, then for $h \in
                \mathbb{R}^d$ satisfying $0<\|h\|< R/2$, we have by
                the Mean Value Theorem that for all $x\in \mathbb{R}^d \setminus B_R(s)$,
		\begin{equation}
                  g(s+h,x) - g(s,x)\; =\; \int_0^1\!\nabla_s g(s+c\,h,x) \,dc \, \cdot \,h. \label{eq: MVT} 
		\end{equation}
		where for $a,b\in\R^d$, $a\cdot b=\langle a,b\rangle$
                denotes the usual vector dot product. Hence,
		\begin{align*}  
			\frac{|g(s+h,x) - g(s,x) - \nabla_s g(s,x) \cdot h|}{\|h\|} & \leq \mynorm{\int_0^1\!\nabla_s g(s+c\,h,x) \,dc} + \| \nabla_s g(s,x)\| \\ 
			&=  \mynorm{-k\,\int_0^1\!\frac{s+c\,h-x}{\|s+c\,h-x\|^{k+2}} \,dc} + k\,\mynorm{ \frac{s-x}{\|s-x\|^{k+2}}} \\ 
			&\leq k\,\int_0^1\!\|s+c\,h-x\|^{-k-1}\,dc + k \,\|s-x\|^{-k-1} .
		\end{align*}
		Integrating over $\mathbb{R}^d \setminus B_R(s)$, we obtain,
		\begin{align}	
			&\int_{\mathbb{R}^d \setminus B_R(s)}\! \frac{|g(s+h,x) - g(s,x) -\nabla_s g(s,x) \cdot h|}{\|h\|} \,\rho(x)\,dV(x)  \\ 
			& \leq \int_{\mathbb{R}^d \setminus B_R(s)}\!k\,\left(\int_0^1\!\|s+c\,h-x\|^{-k-1}\,dc +  \,\|s-x\|^{-k-1}\right)\;\rho(x)\,dV(x)  \\  
			& \leq \int_{\mathbb{R}^d}\!k\,\left(\int_0^1\!\|s+c\,h-x\|^{-k-1}\,dc +  \,\|s-x\|^{-k-1}\right) \;\rho(x)\,dV(x). \label{eq: 2}
		\end{align} 
		By Lemma \ref{Lemma 1}, we have that for all $c\in [0,1]$, 
		\[   \int_{\mathbb{R}^d}\!k\,\|s+c\,h-x\|^{-k-1} \,\rho(x)\,dV(x)  < \infty.  \]
		Defining \[   \theta = \max_{c \in[0,1]}\,\int_{\mathbb{R}^d}\!k\,\|s+c\,h-x\|^{-k-1} \,\rho(x)\,dV(x), \] 
		we see that $\theta$ is a real value (the maximum of a
                continuous function on a compact interval is achieved
                at a point in the interval). Hence 
		we obtain that:
                \[ \int_0^1\int_{\mathbb{R}^d}\!k\,\left(\|s+c\,h-x\|^{-k-1} +  \,\|s-x\|^{-k-1}\right) \,\rho(x)\,dV(x) \,dc \leq  2\, \theta.\]
Now, by \eqref{eq: 2} and Tonelli, note that
		\begin{align*}	 \int_{\mathbb{R}^d \setminus B_R(s)}\! \frac{|g(s+h,x) - g(s,x)-\nabla_s g(s,x) \cdot h|}{\|h\|} \,\rho(x)\,dV(x) 
			& \leq 2\,\theta. 
		\end{align*}
By the Mean Value Theorem applied to \eqref{eq: MVT}, and for each
$s$, we have the following pointwise convergence as $h \rightarrow 0$:
		\begin{align*}  \lim_{\|h\| \rightarrow 0} \frac{|g(s+h,x) - g(s,x) - \nabla_s g(s,x) \cdot h|}{\|h\|} &\leq \lim_{\|h\| \rightarrow 0} \|\nabla_s g(s+c\,h,x)  - \nabla_s g(s,x)\| = 0.
		\end{align*}
		Therefore, by the Dominated Convergence Theorem,
		\begin{align*}
			\lim_{h \rightarrow 0} \, \frac{1}{\|h\|}\,\int_{\mathbb{R}^d \setminus B_R(s)}\! \left(g(s+h,x) - g(s,x)    - \nabla_s g(s,x) \cdot h\right) \,\rho(x)\,dV(x) &= 0.
		\end{align*}
		By the definition of the Frechet derivative, we have:
		\begin{align*} \nabla_s \int_{\mathbb{R}^d \setminus B_R(s) }\! g(s,x)\, \,\rho(x)\,dV(x)   & = \int_{\mathbb{R}^d \setminus B_R(s) }\! \nabla_s g(s,x) \,\rho(x)\,dV(x) \\
			&=  -k\,\int_{\mathbb{R}^d \setminus B_R(s)}\!\frac{(s-x)}{\|s-x\|^{k+2}}\,\rho(x)\,dV(x).\end{align*}
		Noting that this result does not depend on $R > 0$, let $R \rightarrow 0$ to obtain:
		\[
		\lim_{R \rightarrow 0} \,\nabla_s \int_{\mathbb{R}^d \setminus B_R(s)}\!  g(s,x)   \,\rho(x)\,dV(x) =  -k\,\lim_{R \rightarrow 0} \,\int_{\mathbb{R}^d \setminus B_R(s) }\!\frac{(s-x)}{\|s-x\|^{k+2}}\,\rho(x)\,dV(x). \]
 Due to Lemma \ref{Lemma 1}, the integral on the right hand side
 converges absolutely and uniformly as $R \rightarrow 0$, and for any
 $0 \leq k \leq d -2$. Due to the same lemma, the integral on the left
 hand side also converges for each $s$ to $\int_{\mathbb{R}^d }\!
 g(s,x)    \,\rho(x)\,dV(x)$. Therefore, by a standard result in
 analysis (e.g., Theorem 7.17 in \cite{rudin}), the above integrals
 are equal over all of $\R^d$:
		\[
		\nabla_s \int_{\mathbb{R}^d }\!  g(s,x)    \,\rho(x)\,dV(x) =  -k\,\int_{\mathbb{R}^d  }\!\frac{(s-x)}{\|s-x\|^{k+2}}\,\rho(x)\,dV(x).\] 
            \end{proof}
            (A simple modification of the proof also works for integrands the type $\phi(g(s, x))$ where $\phi$ is a differentiable function.)
	
	%%%%%%%%%%%%%%%%%%%%%%%%%%%%%%%%%%%%%%%%%%
	%%%%%%%% Proof of Lemma 3.3
	%%%%%%%%%%%%%%%%%%%%%%%%%%%%%%%%%%%%%%%%%%%
	\section{Proof of Lemma \ref{Lemma 3}}
	\begin{proof}
Once again write $\|\cdot\|$ to mean the usual
                Euclidean norm $\|\cdot\|_2$.		In Lemma~\ref{Lemma 2}, we showed that under the given conditions on $\rho$,  the Frechet derivative of $G$ with respect to $s$ is given by, 
		%%%%%%%%%% Equation 01 %%%%%%%%%%%%%%%%%%%%%%
		\begin{equation}
			\nabla_{s} G(s) = -k \int_{\mathbb{R}^d} \frac{s-x}{\|s-x\|^{k+2}}\rho(x) dV(x). \label{eqn: Grad of F(s)}
		\end{equation}	Below we repeat the steps of the
                proof of Lemma~\ref{Lemma 3}, which hold even in
                the presence of the additional complication that the integrand in \eqref{eqn: Grad of F(s)} is a vector function of $s$.
		Let $R>0$ and
                $g_1(s,x)=-k(s-x)\|s-x\|^{-(k+2)}$, which is bounded and differentiable on
                $\mathbb{R}^d\backslash B_{R}(s)$. Now, for $h \in
                \mathbb{R}^d$ satisfying $0<\|h\|<R/2$, an
                application of the Mean Value Theorem yields that for all $x\in\mathbb{R}^d\backslash B_{R}(s)$, 
		%%%%%%% Equation 03 %%%%%%%%%%%%%%%%%%%%%%
		\begin{equation} 
			g_1(s+h,x)- g_1(s,x)  =\left[\int_{0}^{1} D_{s}g_1\left(s+ch,x\right)dc\right]\cdot h \label{eqn: Mean Value thm}
		\end{equation}
		%%%%%%%%%%%%%%%%%%%%%%%%%%%%%%%%%%%%%%%%
		Noting that,
		\begin{equation*}
			\begin{split}
				& \Bigg \|\left[-k\int_{0}^{1}\left(\frac{I_d}{\|s+ch-x\|^{k+2}}-\frac{(k+2)\left(s+ch-x\right)\left(s+ch-x\right)^{t}}{\|s+ch-x\|^{k+4}}\right)dc\right] \cdot h \Bigg \|\\
				& \leq \int_{0}^{1} k \left(\frac{\sqrt{d}}{\|s+ch-x\|^{k+2}}+(k+2)\frac{\|s+ch-x\|^2}{\|s+ch-x\|^{k+4}}\right)dc\\
				& \leq k(\sqrt{d}+k+2) \int_{0}^{1}\|s+ch-x\|^{-k-2}dc.
			\end{split}
		\end{equation*}
		we obtain, 
		%%%%%%%%%%%%%%%% Equation 04 %%%%%%%%%%%%%%%%
		\begin{align}
			&\frac{|g_1(s+h,x)-g_1(s,x)-D_{s}g_1(s,x).h|}{\|h\|} \leq  \nonumber \\ & \qquad \qquad  k(\sqrt{d}+k+2)\left[\int_{0}^{1}\|s+ch-x\|^{-k-2}dc+\|s-x\|^{-k-2}\right].\label{eqn: Quotient computation}
		\end{align}
		%%%%%%%%%%%%%%%%%%%%%%%%%%%%%%%%%%%%%%%%%%%%%%%%
		Integrating both sides of (\ref{eqn: Quotient computation}) over $\mathbb{R}^d\backslash B_{R}(s)$ yields, 
		%%%%%%%%%%% Equation 05 %%%%%%%%%%%%%%%%%%%%%%
		\begin{equation}
			\begin{split}
				& \int_{\mathbb{R}^d\backslash B_{R}(s)} \frac{|g_1(s+h,x)-g_1(s,x)-D_{s}g_1(s,x).h|}{\|h\|} \rho(x) dV(x)\\
				& \leq k(\sqrt{d}+k+2)\int_{\mathbb{R}^d\backslash B_{R}(s)}\left[\int_{0}^{1}\|s+ch-x\|^{-k-2}dc+\|s-x\|^{-k-2}\right] \rho(x) dV(x)\\
				& \leq k(\sqrt{d}+k+2)\int_{\mathbb{R}^d}\left[\int_{0}^{1}\|s+ch-x\|^{-k-2}dc+\|s-x\|^{-k-2}\right] \rho(x) dV(x).\label{eqn: Integration 1.1}
			\end{split}
		\end{equation}
		%%%%%%%%%%%%%%%%%%%%%%%%%%%%%%%%%%%%%%%
		By Lemma \ref{Lemma 2}, we have that for all $c \in \left[0,1\right]$,  
		%%%%%%%%%%%% Equation 06 %%%%%%%%%%%%
		\begin{equation}
			\int_{\mathbb{R}^d}k(\sqrt{d}+k+2) \|s+ch-x\|^{-k-2}\rho(x)dV(x) < \infty. \label{eqn: Integration 1.2}
		\end{equation}
		%%%%%%%%%%%%%%%%%%%%%%%%%%%%%%%%%%%%%
		Denoting 
		%%%%%%%%%%%%% Equation 07 %%%%%%%%%%
		\begin{equation}
			\nu = \underset{c \in [0,1]}{\max}\int_{\mathbb{R}^d}k(\sqrt{d}+k+2) \|s+ch-x\|^{-k-2}\rho(x)dV(x),\label{eqn: max value}
		\end{equation}
		%%%%%%%%%%%%%%%%%%%%%%%%%%%%%%%%%
		we have immediately that, 
		%%%%%%%%%% Equation 08 %%%%%%%%%
		\begin{equation}
			\int_{0}^{1} \int_{\mathbb{R}^d} k(\sqrt{d}+k+2) \left(\|s+ch-x\|^{-k-2}+\|s-x\|^{-k-2}\right)\rho(x) dV(x) dc \leq 2\nu. \label{eqn: Integration 1.3}
		\end{equation}
		Now, from (\ref{eqn: Integration 1.1}) we have from
                Fubini's Theorem that, 
		\begin{equation}
			\int_{\mathbb{R}^d\backslash B_{R}(s)}\,\frac{|g_1(s+h,x)-g_1(s,x)-D_{s}g_1(s,x) \cdot h|}{\|h\|}\, \rho(x)\,dV(x)\leq 2\nu. \label{eqn: Integration 1.4 }
		\end{equation}
		By \eqref{eqn: Mean Value thm}, for each $s$, the function $( g_1(s+h,x) - g_1(s,x) - D_s g_1(s,x) \cdot h)/\|h\|$ converges pointwise to $0$ as $h \rightarrow 0$, since
		\begin{align*}  
			\lim_{\|h\| \rightarrow 0} \frac{|g_1(s+h,x) - g_1(s,x) - D_s g_1(s,x) \cdot h|}{\|h\|} &\leq \lim_{\|h\| \rightarrow 0} \Bigg \|\int_0^1\!D_s g_1(s+c\,h,x) \,dc  - D_s g_1(s,x) \Bigg \| \\ 
			&=  \Bigg \| \int_0^1\!D_s g_1(s,x) \,dc  - D_s g_1(s,x) \Bigg \|\\
			&= 0,
		\end{align*}
		whence, by Lebesgue's Dominated Convergence Theorem, 
		\begin{equation}
			\underset{h \rightarrow 0}{\lim} \int_{\mathbb{R}^d\backslash B_{R}(s)} \frac{|g_1(s+h,x)-g_1(s,x)-D_{s}g_1(s,x) \cdot h|}{\|h\|}\rho(x)dV(x)=0.\label{eqn: Integration 1.5}
		\end{equation}
		By the definition of the Frechet derivative, 
		\begin{align*}
			& D_{s} \int_{\mathbb{R}^d\backslash B_{R}(s)} g_1(s,x) \, \rho(x) \, dV(x) 
			=  \int_{\mathbb{R}^d\backslash B_{R}(s)}  \!D_{s}g_1(s,x)  \,\rho(x)\, dV(x) \\
			&=  -k \int_{\mathbb{R}^d\backslash B_{R}(s)} \!\!\! \left(\frac{I_d}{\|s-x\|^{k+2}}-(k+2)\frac{(s-x)(s-x)^{t}}{\|s-x\|^{k+4}}\right)\rho(x)dV(x). \label{eqn: Integration 1.6}
		\end{align*}
	Finally, letting $R \rightarrow 0$ yields the desired result
        by  the same argument used at the end of Lemma \ref{Lemma 2}.
	\end{proof}
	%%%%%%%%%%%%%%%%%%%%%%%%%%%%%%%%%%%%%%%%%%
	%%%%%%%% Proof of Lemma 3.4
	
	%%%%%%%%%%%%%%%%%%%%%%%%%%%%%%%%%%%%%%%%%%%
	\section{Proof of Lemma \ref{Lemma 4}}
		\begin{proof}%[Proof of Lemma \ref{Lemma 4}]
		For $0 < \varepsilon < 1 $, there exists $[t_1,t_2]$ such that $\mu([t_1,t_2]) > 1 - \varepsilon.$ 
		Let $s \in [t_1,t_2]$ and $\{s_k\}_{k=1}^{\infty}$ be a sequence converging to $0$. Define, 
		\begin{align}
			f_{\mu}(s) &= \int_{\mathbb{R}} (|s-x|-|x| )\rho(x) dx,  \\
			\phi(s;t_1,t_2) &= \int_{[t_1,t_2]} (|s-x|-|x| )\rho(x) dx. 
		\end{align}
		As $\rho$ is a pdf, for each $s$ the function $\theta(x) := (|s-x| - |x|)\,\rho(x)$ $\in L^1([t_1,t_2])$. Therefore,
\[
			\frac{\phi(s+s_k;t_1,t_2)-\phi(s;t_1,t_2)}{s_k}
                        = \frac{1}{s_k}\, \int_{[t_1,t_2]}
                        (|s+s_k-x|-|s-x|)\rho(x) dx := A+B,
                      \]
where
\[
A= \frac{1}{s_k}\, \left[ \int_{t_1}^{s+s_k} \!(s+s_k-x)\,\rho(x)\,dx
  -   \int_{t_1}^{s} \!(s-x)\,\rho(x)\,dx\right],
\]
and
\[
B= \frac{1}{s_k}\, \left[ \int_{s+s_k}^{t_2} \!(s+s_k-x)\,\rho(x)\,dx -   \int_{s}^{t_2} \!(s-x)\,\rho(x)\,dx\right]. 
\]
Now note that		
		\begin{align*}  
A &=  \frac{1}{s_k} \,\int_{t_1}^{s}\,s_k\rho(x)dx + \frac{1}{s_k} \int_{s}^{s+s_k} \!\!\!\!\! (s+s_k-x)\,\rho(x)dx \\
			&  = \mu([t_1,s]) + \frac{1}{s_k} \int_{s}^{s+s_k} \!\!\!\!\! (s+s_k-x)\,\rho(x)\,dx ,	
		\end{align*}
		and similarly,
		\[ B =  -\mu([s,t_2]) + \frac{1}{s_k} \int_{s}^{s+s_k}
                  \!\!\!\!\! (s+s_k-x)\,\rho(x)\,dx.  \]
                Putting things together,
		\[  \frac{\phi(s+s_k;t_1,t_2)-\phi(s;t_1,t_2)}{s_k} = \mu([t_1,s]) -\mu([s,t_2]) + \frac{2}{s_k} \int_{s}^{s+s_k} \!\!\!\!\! (s+s_k-x)\,\rho(x)\,dx.\] 
		By Theorem~\ref{th:Leibniz rule} applied to $(s,y) \in
                [t_1,t_2] \times [t_1,t_2]$, yields that, for almost all $s \in [t_1,t_2]$, 
		\begin{equation}\frac{d}{dy} \Big{\vert}_{y = s}\,\int_{s}^{y} \! (y-x)\,\rho(x)\,dx = 0 \;\;\;\Longrightarrow\;\; \lim_{s_k \rightarrow 0}\,\frac{1}{s_k} \int_{s}^{s+s_k} \! (s+s_k-x)\,\rho(x)\,dx = 0. \label{eq:leibiz} \end{equation} It follows that for almost all $s \in [t_1,t_2]$,
		\[ \phi'(s) = \mu([t_1,s]) -\mu([s,t_2]) = \mu([-\infty,s]) -\mu([s,\infty]) + 2\,\varepsilon.\] Letting $\varepsilon \rightarrow 0$, we see that  for almost all $s \in \mathbb{R}$,
		$f'_\mu(s) = \mu([-\infty,s]) -\mu([s,\infty])$.
                Now, applying Leibnitz's rule (2nd fundamental theorem of calculus), yields that, for almost all $s \in \mathbb{R}$,
		\begin{equation} f''_\mu(s) = \rho(s) - (-\rho(s)) =
                  2\,\rho(s), \label{eq:lemma3_4
                    conclusion}\end{equation}
                thus proving the desired result. If $\rho$ is a continuous function, then (the
                same Riemann integral version of) Leibnitz's rule may be employed in \eqref{eq:leibiz} instead of Theorem \ref{th:Leibniz rule} to conclude \eqref{eq:lemma3_4 conclusion} for all $s \in \mathbb{R}$.
	\end{proof}
	
	\section{Proof of Theorem \ref{Theorem}}
		By Lemma \ref{Lemma 2}, we can obtain the gradient of $f_{\mu}(s)$ as, 
		\begin{equation}
			\nabla_{s}f_{\mu}(s) = \E\left( \frac{s-X}{\|s-X\|_2} \right), \label{eqn: Proof 2.1}
		\end{equation}
		and by Lemma \ref{Lemma 3} we can compute its Hessian as
		\begin{equation}
			D_s\,\nabla_{s} f_{\mu}(s) = \E\left(\frac{I_d}{\|s-X\|_2}-\frac{(s-X)(s-X)^t}{\|s-X\|_2^3}\right), \label{eqn: Proof 2.2}
		\end{equation}
		which we note is symmetric and positive definite. Hence, we can compute the Laplacian of $f_{\mu}$ as
		\begin{equation}
				\Delta_{s} f_{\mu}(s) = \text{trace}\left(D_s\,\nabla_{s} f_{\mu}(s)\right) =(d-1) \E\left(\frac{1}{\|s-X\|_2}\right).\label{eqn: Proof 2.3}
		\end{equation}
		Now, by Lemma \ref{Lemma 1}, $\Delta_{s} f_{\mu}(s)$
                exists for $d \geq 3$, so that \eqref{eqn: Theorem}
                holds for $j = 1$. Proceeding by induction,
		we suppose that \eqref{eqn: Theorem} holds for all
                $j\,\in \{1,\ldots,N-1 \}$, and aim to show that
                (\ref{eqn: Theorem}) holds also for $j+1$.
                
		To this effect, set $k = 2\,j - 1$, and using Lemmas
                \ref{Lemma 1} and \ref{Lemma 3}, we compute, respectively, the
                Frechet gradient and Hessian of
                $\E\|s-X\|_2^{-k}$ as: 
		\begin{equation}
			\nabla_{s} \E\left(\frac{1}{\|s-X\|_2^k}\right) = -k\, \E\left(\frac{s-X}{\|s-X\|_2^{k+2}}\right), \label{eqn: Proof 2.4}
		\end{equation}
		and 
		\begin{equation}
			D_s\nabla_{s}\E\left(\frac{1}{\|s-X\|_2^k}\right) = -k\, \E\left(\frac{I_{d}}{\|s-X\|_2^{k+2}}-(k+2)\frac{(s-X)(s-X)^{t}}{\|s-X\|_2^{k+4}}\right). \label{eqn: Proof 2.5}
		\end{equation}
		Hence, the Laplacian of $\E\|s-X\|_2^{-k}$ is given by:
		\begin{equation}
			\Delta_{s} \E\left(\frac{1}{\|s-X\|_2^k}\right) = -k\,(d-k-2)\, \E\left(\frac{1}{\|s-X\|_2^{k+2}} \right). \label{eqn: Proof 2.6}
                      \end{equation}
                      Now, the expectation on the right hand side of
                      \eqref{eqn: Proof 2.6} exists by Lemma
                      \ref{Lemma 2}, as long as $k \, \leq \, d-2$, or
                      equivalently, $j\leq N$.
		Converting back to the $j$ variable, as
                $\Delta^{j}_{s}f_{\mu}(s)  =
                \Delta_{s}\,(\Delta^{j-1}_{s}f_{\mu}(s))$, yields the
                final result.
	
	%%%%%%%%%%%%%%%%%%%%%%%%%%%%%%%%%%%%%
	%%%%%%%      Corollary 
	%%%%%%%%%%%%%%%%%%%%%%%%%%%%%%%%%%%%
	\section{Proof of Corollary \ref{corollary}}
	\begin{proof}
		Consider the Poisson equation, 
		\begin{equation}
			\Delta \phi(s) = \rho(s), \quad \text{a.e.}, \label{eqn: Poisson 1.1}
		\end{equation}
		whose fundamental solution is
		\[  
		\phi(s) = \frac{\Gamma\left(\frac{d}{2}+1\right)}{d\,(d-2)\pi^{d/2}}\E\left( \frac{1}{\|s-X\|_2^{d-2}} \right).
		\]
		Putting \eqref{eq: final polyharmonic} and \eqref{eqn:
                  Poisson 1.1} together, we obtain:
		\begin{align*}
			\Delta^{N} f_{\mu}(s) &= (-1)^{N-1}\,(d-1)\,\cdots\,(2)\E\left(\frac{1}{\|s-X\|_2^{d-2}} \right) \\
			&= (-1)^{N-1}\,(d-1)\,\cdots\,(2)\, d\,(d-2)\,\frac{\pi^{d/2}}{\Gamma\left(\frac{d}{2}+1\right)}\,\phi(s),
		\end{align*}
                so that
		\begin{align*}
			\Delta\,\Delta^{N} f_{\mu}(s) &=  (-1)^{N-1}\,(d-1)\,\cdots\,(2)\E\left(\frac{1}{\|s-X\|_2^{d-2}} \right) \\
			&= (-1)^{N-1}\,(d-1)\,\cdots\,(2)\, d\,(d-2)\,\frac{\pi^{d/2}}{\Gamma\left(\frac{d}{2}+1\right)}\, \rho(s), \quad \text{a.e},
		\end{align*}
                thus proving the result.
	\end{proof}
	
\end{appendix}
\bibliographystyle{imsart-number} % Style BST file (imsart-number.bst or imsart-nameyear.bst)
\bibliography{references.bbl}       % Bibliography file (usually '*.bib')

%% or include bibliography directly:
% \begin{thebibliography}{}
% \bibitem{b1}
% \end{thebibliography}

\end{document}